\documentclass[a4paper]{article}
\usepackage{amsfonts}
\usepackage{amssymb}
\usepackage{amsmath}
\usepackage{amscd}
\usepackage{latexsym}
\usepackage{epsfig}
\usepackage{calrsfs}
\newtheorem{theorem}{Theorem}[section]
\newtheorem{conjecture}[theorem]{Conjecture}
\newtheorem{corollary}[theorem]{Corollary}
\newtheorem{lemma}[theorem]{Lemma}
\newtheorem{proposition}[theorem]{Proposition}
\newtheorem{question}[theorem]{Question}
\newtheorem{example}[theorem]{Example}
\newtheorem{definition}[theorem]{Definition}
\newtheorem{remark}[theorem]{Remark}
\newtheorem{problem}[theorem]{Problem}

\newcommand{\proof}{\medskip \noindent {\bf Proof. \ \ }}
\newcommand{\qed}{\null\hfill $\Box\;\;$ \medskip}
\begin{document}

\parbox{1mm}

\begin{center}
{\bf {\sc \Large Flett's mean value theorem: a survey}}
\end{center}

\vskip 12pt

\begin{center}
{\bf Ondrej HUTN\'IK\footnote{corresponding author} and Jana MOLN\'AROV\'A} \\
Institute of Mathematics, Faculty of Science, \\ Pavol Jozef \v
Saf\'arik University in Ko\v sice, \\ Jesenn\'a 5, SK 040~01 Ko\v
sice, Slovakia,
\\ {\it E-mail address:} ondrej.hutnik@upjs.sk, jana.molnarova88@gmail.com
\end{center}

\vskip 24pt

\hspace{5mm}\parbox[t]{12cm}{\fontsize{9pt}{0.1in}\selectfont\noindent{\bf
Abstract.} This paper reviews the current state of the art of the
mean value theorem due to Thomas M. Flett. We present the results
with detailed proofs and provide many new proofs of known results.
Moreover, some new observations and yet unpublished results are
included.} \vskip 12pt

\hspace{5mm}\parbox[t]{12cm}{\fontsize{9pt}{0.1in}\selectfont\noindent{\bf
Key words and phrases.} Flett's mean value theorem, real-valued
function, differential function, Taylor's polynomial,
Pawlikowska's theorem} \vskip 12pt

\hspace{5mm}\parbox[t]{12cm}{\fontsize{9pt}{0.1in}\selectfont\noindent{\bf
Mathematics Subject Classification (2010)} 26A24, 26D15}

\vskip 24pt

%%%%%%%%%%%%%%%%%%%%%%%%%%%%%%%%%%%%%%%%%%%%%%%%%%%%%%%%%%%%%%%%%%%%%%%%%%%%%%%%%%%%%%%%%%%%%%%
\section{Introduction and preliminaries}
%%%%%%%%%%%%%%%%%%%%%%%%%%%%%%%%%%%%%%%%%%%%%%%%%%%%%%%%%%%%%%%%%%%%%%%%%%%%%%%%%%%%%%%%%%%%%%%

\begin{figure}
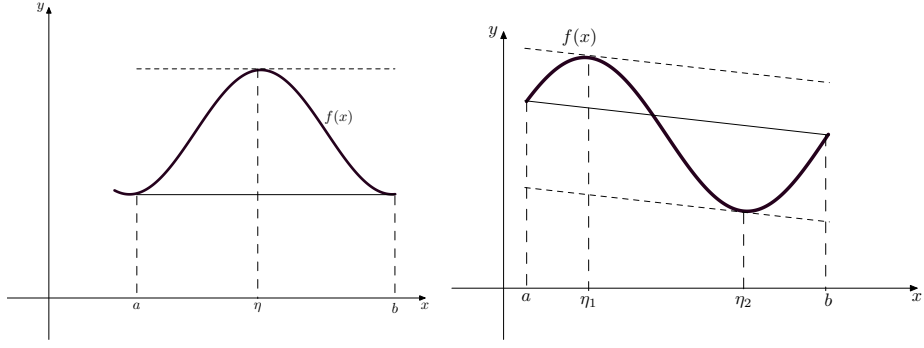
%[hbp!]
\begin{center}
\includegraphics[scale=0.55]{Rolleovaveta.1}\hskip 8pt
\includegraphics[scale=0.68]{Lagrangeovaveta.1}
\caption{Geometrical interpretation of Rolle's and Lagrange's
theorem} \label{obrlagrange}
\end{center}
\end{figure}

\paragraph{Motivations and basic aim} Mean value theorems of differential and integral calculus
provide a relatively simple, but very powerful tool of
mathematical analysis suitable for solving many diverse problems.
Every student of mathematics knows the Lagrange's mean value
theorem which has appeared in Lagrange's book \emph{Theorie des
functions analytiques} in 1797 as an extension of Rolle's result
from 1691. More precisely, Lagrange's theorem says that for a
continuous (real-valued) function $f$ on a compact set $\langle
a,b\rangle$ which is differentiable on $(a,b)$ there exists a
point $\eta\in(a,b)$ such that
$$f'(\eta)=\frac{f(b)-f(a)}{b-a}.$$ Geometrically Lagrange's theorem states that
given a line $\ell$ joining two points on the graph of a
differentiable function $f$, namely $[a, f(a)]$ and $[b, f(b)]$,
then there exists a point $\eta\in (a, b)$ such that the tangent
at $[\eta, f(\eta)]$ is parallel to the given line $\ell$, see
Fig.~\ref{obrlagrange}. Clearly, Lagrange's theorem reduces to
Rolle's theorem if $f(a)=f(b)$. In connection with these
well-known facts the following questions may arise: \textsl{Are
there changes if in Rolle's theorem the hypothesis $f(a) = f(b)$
refers to higher-order derivatives? Then, is there any analogy
with the Lagrange's theorem? Which geometrical consequences do
such results have?} These (and many other) questions will be
investigated in this paper in which we provide a survey of known
results as well as of our observations and obtained new results.

\paragraph{Notation} Throughout this paper we will use the following unified notation:
$\mathcal{C}(M)$, resp. $\mathcal{D}^n(M)$, will denote the spaces
of continuous, resp. $n$-times differentiable real-valued
functions on a set $M\subseteq \mathbb{R}$. Usually we will work
with a compact set of the real line, i.e., $M=\langle a,b\rangle$
with $a<b$. Therefore, we recall that under continuity of a
function on $\langle a,b\rangle$ we understand its continuity on
$(a,b)$ and one-sided continuity at the end points of the
interval. Similarly we will understand the notion of
differentiability on a closed interval. For functions $f,g$ on an
interval $\langle a,b\rangle$ (for which the following expression
has its sense) the expressions of the form
$$\frac{f^{(n)}(b)-f^{(n)}(a)}{g^{(n)}(b)-g^{(n)}(a)}, \quad
n\in\mathbb{N}\cup\{0\},$$ will be denoted by the symbol
$\phantom{}_a^b\mathcal{K}\left(f^{(n)},g^{(n)}\right)$. If the
denominator is equal to $b-a$, we will write only
$\phantom{}_a^b\mathcal{K}\left(f^{(n)}\right)$. So, Lagrange's
theorem in the introduced notation has the following form: if
$f\in\mathcal{C}\langle a,b\rangle\cap\mathcal{D}(a,b)$, then
there exists $\eta\in(a,b)$ such that
$f'(\eta)=\phantom{}_a^b\mathcal{K}(f)$, where we use the usual
convention $f^{(0)}:=f$.

\paragraph{Structure of this paper} The organization of this paper
is as follows: in Section~\ref{secFlett} we present the original
result of Flett as well as its generalization due to Riedel and
Sahoo removing the boundary condition. Further sufficient
conditions of Trahan and Tong for validity of assertion of Flett's
theorem are described in Section~\ref{sectionsufficient} together
with proving two new extensions and the detailed comparison of all
the presented conditions. Section~\ref{secintegral} deals with
integral version of Flett's theorem and related results. In the
last Section~\ref{secPawlikowska} we give a new proof of
higher-order generalization of Flett's mean value theorem due to
Pawlikowska and we present a version of Flett's and Pawlikowska's
theorem for divided differences of a real function.

\section{Flett's mean value theorem}\label{secFlett}

Let us begin with the following easy observation
from~\cite{Flett}: if $g\in\mathcal{C}\langle a,b\rangle$, then
from the integral mean value theorem there exists $\eta\in(a,b)$
such that
$$g(\eta)=\frac{1}{b-a}\int_{a}^{b} g(t)\,\mathrm{d}t.$$ Moreover, if we consider the function $g\in\mathcal{C}\langle
a,b\rangle$ with the properties
\begin{equation}\label{podmienka}
g(a)=0, \quad \int_{a}^{b}g(t)\,\mathrm{d}t=0,
\end{equation} and define the function
$$\varphi(x) =
\begin{cases}
\frac{1}{x-a}\int_{a}^{x}g(t)\,\mathrm{d}t, & x\in(a,b\rangle, \\
\nonumber 0, & x=a,
\end{cases}$$
then $\varphi\in\mathcal{C}\langle a,b\rangle\cap\mathcal{D}(a,b)$
and $\varphi(a)=0=\varphi(b)$. Thus, by Rolle's theorem there
exists $\eta\in(a,b)$ such that $\varphi'(\eta)=0$, i.e.,
$$\frac{g(\eta)}{\eta-a}-\frac{1}{(\eta-a)^2} \int_{a}^{\eta} g(t)\,\mathrm{d}t = 0\,\,\,\Leftrightarrow
\,\,\, g(\eta)=\frac{1}{\eta-a}\int_{a}^{\eta}
g(t)\,\mathrm{d}t.$$ The latter formula resembles the one from
integral mean value theorem replacing formally $b$ by $\eta$. It
is well-known that the mean value $\varphi$ (known as the integral
mean) of function $g$ on the interval $\langle a,x\rangle$ is in
general less irregular in its behaviour than $g$ itself. When
defining the function $g$ we may ask whether the second condition
in~(\ref{podmienka}) may be replaced by a simpler condition, e.g.,
by the condition $g(b)=0$. Later we will show that it is possible
and the result in this more general form is a consequence of
Darboux's intermediate value theorem, see the first proof of
Flett's theorem. If we define the function
$$f(x)=\int_{a}^{x}g(t)\,\mathrm{d}t, \quad x\in\langle a,b\rangle,$$
from our considerations we get an equivalent form of result to
which this paper is devoted. This result is an observation of
\textsc{Thomas Muirhead Flett} (1923--1976) from 1958 published in
his paper~\cite{Flett}. Indeed, it is a variation on the theme of
Rolle's theorem where the condition $f(a)=f(b)$ is replaced by
$f'(a)=f'(b)$, or, we may say that it is a Lagrange's type mean
value theorem with a Rolle's type condition.

\begin{theorem}[Flett, 1958]\label{flettova veta}
If $f\in\mathcal{D}\langle a,b\rangle$ and $f'(a)=f'(b)$, then
there exists $\eta\in(a,b)$ such that
\begin{equation}\label{flettMV}
f'(\eta)=\phantom{}_a^\eta\mathcal{K}(f).
\end{equation}
\end{theorem}

For the sake of completeness we give the original proof of Flett's
theorem adapted from~\cite{Flett} and rewritten in the sense of
introduced notation.

\vskip 10pt\noindent \textbf{Proof of Flett's theorem I.}\,\,\,\,
Without loss of generality assume that $f'(a)=f'(b)=0$. If it is
not the case we take the function $h(x)=f(x)-xf'(a)$ for
$x\in\langle a,b\rangle$. Put
\begin{equation}\label{flettfunkcia}
g(x) =
\begin{cases}
\phantom{}_a^x\mathcal{K}(f), & x\in(a,b\rangle \\
f'(a), & x=a.
\end{cases}
\end{equation}
Obviously, $g\in\mathcal{C}\langle a,b\rangle\cap\mathcal{D}(a,b)$
and
$$g'(x)=-\frac{\phantom{}_a^x\mathcal{K}(f)-f'(x)}{x-a} = -\phantom{}_a^x\mathcal{K}(g)+\phantom{}_a^x\mathcal{K}(f'),
\quad x\in(a,b\rangle.$$ It is enough to show that there exists
$\eta\in(a,b)$ such that $g'(\eta)=0$.

From the definition of $g$ we have that $g(a)=0$. If $g(b)=0$,
then Rolle's theorem guarantees the existence of a point
$\eta\in(a,b)$ such that $g'(\eta)=0$. Let $g(b)\neq 0$ and
suppose that $g(b)>0$ (similar arguments apply if $g(b)<0$). Then
$$g'(b)=-\phantom{}_a^b\mathcal{K}(g)=-\frac{g(b)}{b-a}<0.$$
Since $g\in\mathcal{C}\langle a,b\rangle$ and $g'(b)<0$, i.e., $g$
is strictly decreasing in $b$, then there exists $x_{1}\in(a,b)$
such that $g(x_{1})>g(b)$. From continuity of $g$ on $\langle
a,x_{1}\rangle$ and from relations $0=g(a)<g(b)<g(x_{1})$ we
deduce from Darboux's intermediate value theorem that there exists
$x_{2}\in(a,x_{1})$ such that $g(x_{2})=g(b)$. Since
$g\in\mathcal{C}\langle x_2,b\rangle\cap\mathcal{D}(x_2,b)$, from
Rolle's theorem we have $g'(\eta)=0$ for some
$\eta\in(x_2,b)\subset(a,b)$. \qed

A different proof of Flett's theorem using Fermat's theorem
(necessary condition for the existence of a local extremum) may be
found in~\cite[p.225]{RRA}.

\vskip 10pt\noindent \textbf{Proof of Flett's theorem II.}\,\,\,\,
Let us consider the function $g$ defined by~(\ref{flettfunkcia}).
If $g$ achieve an extremum at an interior point $\eta\in(a,b)$,
then Fermat's theorem yields $g'(\eta)=0$ and we conclude the
proof.

Assume the contrary, i.e., $g$ achieves an extremum only at the
point $a$ or $b$. Without loss of generality we may assume that
for each $x\in\langle a,b\rangle$ we have $g(a)\leq g(x)\leq
g(b)$. From the second inequality we get
$$(\forall x\in\langle a,b\rangle)\,\, f(x)\leq f(a)+(x-a)g(b).$$ It follows that for each $x\in\langle a,b)$
we have
$$\phantom{}_x^b\mathcal{K}(f)\geq \frac{f(b)-f(a)-(x-a)g(b)}{b-x}=\phantom{}_a^b\mathcal{K}(f).$$
If $x\rightarrow b^-$, then $f'(b)\geq g(b)$ which yields
$f'(a)\geq g(b)$. But $f'(a)=g(a)$, so $g(a)\geq g(b)$. This
implies that $g$ is constant on $(a,b)$, that is $g'(x)=0$ for
each $x\in(a,b)$. Then for each $\eta\in(a,b)$ we have
$f'(\eta)=\phantom{}_a^\eta\mathcal{K}(f)$, which finishes the
proof. \qed

\begin{figure}
\begin{center}
\includegraphics[scale=0.8]{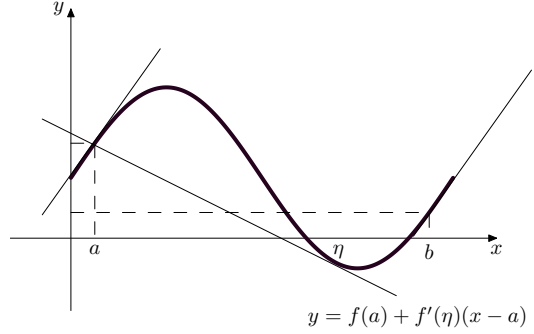}
\caption{Geometrical interpretation of Flett's theorem}
\label{obrFlett1}
\end{center}
\end{figure}

\paragraph{Geometrical meaning of Flett's theorem} If a
curve $y=f(x)$ has a tangent at each point of $\langle a,b\rangle$
and tangents at the end points $[a,f(a)]$ and $[b,f(b)]$ are
parallel, then Flett's theorem guarantees the existence of such a
point $\eta\in(a,b)$ that the tangent constructed to the graph of
$f$ at that point passes through the point $[a, f(a)]$, see
Fig.~\ref{obrFlett1}.

\begin{example}\rm
In which point of the curve $y=x^3$ the tangent passes through the
point $X=[-2,-8]$?

It is easy to verify that $X$ lies on the curve and $y=x^3$ is
differentiable on $\mathbb{R}$. Since its derivative $y'=3x^{2}$
is even function on $\mathbb{R}$, consider such interval $\langle
a,b\rangle$ to be able to apply Flett's theorem, e.g. $\langle
-2,2\rangle$. Then there exists point (or, points) $\eta\in(-2,2)$
such that
$$3\eta^{2}(\eta+2)=\eta^{3}-(-2)^{3}\,\,\, \Leftrightarrow\,\,\,
\eta^{2}+3\eta-4=0\,\,\, \Leftrightarrow\,\,\,
(\eta+4)(\eta-1)=0.$$ Because $-4\notin(-2,2)$, we consider only
$\eta=1$. Then $y(\eta)=1$ and the desired point is $T=[1,1]$.
\end{example}

\begin{remark}\rm
Clearly, the assertion of Flett's theorem may be valid also in the
cases when its assumption is not fulfilled. For instance, function
$f(x)=|x|$ on the interval $\langle a,b\rangle$, with $a<0<b$, is
not differentiable, but there exist infinite many points
$\eta\in(a,0)$ for which the tangent constructed in the point
$\eta$ passes through the point $[a,-a]$ (since the tangent
coincides with the graph of function $f(x)$ for $x\in(a,0)$).

Another example is the function $g(x)=\textrm{sgn}\,x$ and
$h(x)=[x]$ (sign function and floor function) on the interval
$\langle -1,1\rangle$ which are not differentiable on $\langle
-1,1\rangle$. Finally, the function $k(x)=\arcsin x$ on $\langle
-1,1\rangle$ is not differentiable at the end points, but
assertion of Flett's theorem still holds (we will consider other
sufficient conditions for validity of~(\ref{flettMV}) in
Section~\ref{sectionsufficient}, namely $k$ fulfills Tong's
condition).

We can observe that the functions $g$ and $k$ have improper
derivatives at the points in which are not differentiable, i.e.,
$g'_+(0)=g'_-(0)=k'_+(-1)=k'_-(1)=+\infty$. Therefore, we state
the conjecture that Flett's theorem still holds in that case.

\begin{conjecture}
If $f$ has a proper or improper derivative at each point of the
interval $\langle a,b\rangle$ and the tangents at the end points
are parallel, then there exists $\eta\in(a,b)$ such
that~(\ref{flettMV}) holds.
\end{conjecture}
\end{remark}

\begin{figure}
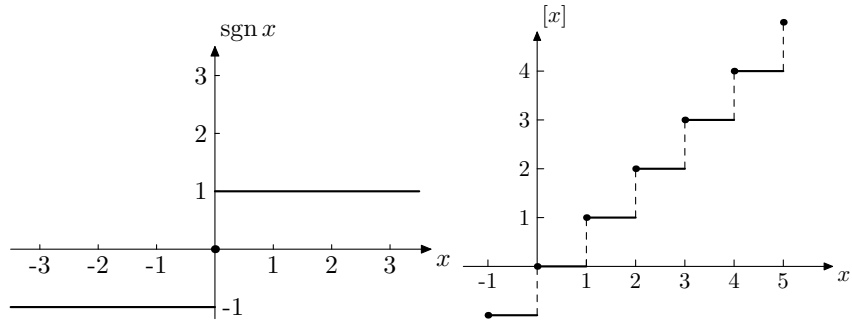

\begin{center}
\includegraphics[scale=0.96]{abssgn.1}
\includegraphics[scale=0.85]{celacast.1}
\caption{Non-differentiable functions for which assertion of
Flett's theorem is valid on the interval $\langle -1,1\rangle$}
\end{center}
\end{figure}

\begin{remark}\rm
Assertion of Flett's theorem may be written in the following
equivalent forms:
$$f'(\eta)= \frac{f(\eta)-f(a)}{\eta-a}
\,\,\,\, \Leftrightarrow \,\,\,\, f(a) =
T_1(f,\eta)(a) \,\,\,\, \Leftrightarrow \,\,\,\, \left|\begin{array}{ccc} f'(\eta) & 1 & 0 \\ f(a) & a & 1 \\
f(\eta) & \eta & 1
\end{array}\right| = 0.$$ In the second expression $T_1(f,x_0)(x)$ is the first
Taylor's polynomial (or, in other words a tangent) of function $f$
at the point $x_0$ as a function of $x$. The last expression
resembles an equivalent formulation of the assertion of Lagrange's
theorem in the form of determinant, i.e.,
$$\left|\begin{array}{ccc}
f'(\eta) & 1 & 0 \\ f(a) & a & 1 \\
f(b) & b & 1
\end{array}\right| = 0.$$ This motivates us to state the following question:
\end{remark}

\begin{question}
Is it possible to find a similar proof (as a derivative of a
function given in the form of determinant) of Flett's theorem?
\end{question}

In connection with applicability of Flett's theorem there exists
many interesting problems proposed and solved by various authors,
see e.g. the problems and solutions section of journals as
American Mathematical Monthly, Electronic Journal of Differential
Equations, etc. A nice application of Flett's theorem for
investigating some integral mean value theorems is given
in~\cite{LL} and similar approach is used in~\cite{GL}. We give
here only one representative example of this kind. The problem
(2011-4 in Electronic Journal of Differential Equations) was
proposed by Duong Viet Thong, Vietnam. The solution to this
problem is our own.

\begin{problem}
Let $f\in\mathcal{C}\langle 0,1\rangle$ and
$$\int_{0}^1 f(x)\,\mathrm{d}x = \int_{0}^1 xf(x)\,\mathrm{d}x.$$ Prove that there exists
$\eta\in(0,1)$ such that
$$\eta^2 f(\eta)=\int_{0}^{\eta} xf(x)\,\mathrm{d}x.$$
\end{problem}

\vskip 10pt\noindent \textbf{Solution.}\,\,\,\, Consider the
differentiable function
$$G(t)=\int_{0}^t xf(x)\,\mathrm{d}x, \quad t\in\langle 0,1\rangle.$$ Clearly, $G'(t)=tf(t)$ for each $t\in\langle 0,1\rangle$.
By~\cite[Lemma 2.8]{LL} there exists $\zeta\in(0,1)$ such that
$G(\zeta)=\int_{0}^\zeta xf(x)\,\mathrm{d}x=0$. Since $G(0)=0$,
then by Rolle's theorem there exists $\theta\in(0,\zeta)$ such
that $G'(\theta)=0$. From $G'(0)=0=G'(\theta)$ by Flett's theorem
there exists $\eta\in(0,\theta)$ such that
$$\hspace{1.6cm} G'(\eta)=\phantom{}_0^{\eta}\mathcal{K}(G)\,\,\, \Leftrightarrow \,\,\, \eta f(\eta) =
\frac{G(\eta)}{\eta}\,\,\, \Leftrightarrow \,\,\, \eta^2
f(\eta)=\int_{0}^{\eta} xf(x)\,\mathrm{d}x.\hspace{1.6cm}\square$$

Naturally, we may ask whether the Lagrange's idea to remove the
equality $f(a)=f(b)$ from Rolle's theorem is applicable for
Flett's theorem, i.e., whether the assumption $f'(a)=f'(b)$ may be
removed for the purpose to obtain a more general result. First
result of that kind has appeared in the book~\cite{RS}.

\begin{theorem}[Riedel-Sahoo, 1998]\label{thmRS}
If $f\in\mathcal{D}\langle a,b\rangle$, then there exists
$\eta\in(a,b)$ such that
$$\phantom{}_a^\eta\mathcal{K}(f)=f'(\eta)-\phantom{}_a^b\mathcal{K}(f')\cdot\frac{\eta-a}{2}.$$
\end{theorem}

In their original proof~\cite{RS} Riedel and Sahoo consider the
auxiliary function $\psi$ given by
\begin{equation}\label{sahoofcia}
\psi(x)=f(x)-\phantom{}_a^b\mathcal{K}(f')\cdot\frac{(x-a)^2}{2},
\quad x\in\langle a,b\rangle,
\end{equation} and apply Flett's theorem to it. Indeed, function $\psi$ is constructed as a difference of $f$
and its quadratic approximation $A+B(x-a)+C(x-a)^2$ at a
neighbourhood of $a$. From $\psi'(a)=\psi'(b)$ we get
$C=\frac{1}{2}\cdot\phantom{}_a^b\mathcal{K}(f')$, and because $A$
and $B$ may be arbitrary, they put $A=B=0$. Of course, the
function $\psi$ is not the only function which does this job. For
instance, the function
$$\Psi(x)=f(x)-\phantom{}_a^b\mathcal{K}(f')\cdot\left(\frac{x^2}{2}-ax\right),
\quad x\in\langle a,b\rangle,$$ does the same, because
$\Psi'(x)=\psi'(x)$ for each $x\in\langle a,b\rangle$. In what
follows we provide a different proof of Riedel-Sahoo's theorem
with an auxiliary function of different form.

\vskip 10pt\noindent \textbf{New proof of Riedel-Sahoo's
theorem.}\,\,\,\, Let us consider the function $F$ defined by
$$F(x)=\left|\begin{array}{cccc}
f(x) & x^2 & x & 1 \\ f(a) & a^2 & a & 1 \\ f'(a) & 2a & 1 & 0 \\
f'(b) & 2b & 1 & 0
\end{array}\right|,\quad x\in\langle a,b\rangle.$$ Clearly,
$F\in\mathcal{D}\langle a,b\rangle$ and
$$F'(x)=\left|\begin{array}{cccc}
f'(x) & 2x & 1 & 0 \\ f(a) & a^2 & a & 1 \\ f'(a) & 2a & 1 & 0 \\
f'(b) & 2b & 1 & 0
\end{array}\right|,\quad x\in\langle a,b\rangle.$$ Thus,
$F'(a)=F'(b)=0$, and by Flett's theorem there exists
$\eta\in(a,b)$ such that
$F'(\eta)=\phantom{}_a^\eta\mathcal{K}(F)$, which is equivalent to
the assertion of Riedel-Sahoo's theorem.\qed

\begin{remark}\rm
As in the case of Flett's theorem it is easy to observe that the
assertion of Riedel-Sahoo's theorem may be equivalently written as
follows $$\left|\begin{array}{ccc} f'(\eta) & 1 & 0 \\ f(a) & a & 1 \\
f(\eta) & \eta & 1
\end{array}\right| =
\phantom{}_a^b\mathcal{K}(f')\cdot\frac{(\eta-a)}{2}.$$
\end{remark}

The geometrical fact behind Flett's theorem is a source of
interesting study in~\cite{DV} we would like to mention here in
connection with Riedel-Sahoo's theorem. Following~\cite{DV} we
will say that \textit{the graph of $f\in\mathcal{C}\langle
a,b\rangle$ intersects its chord in the extended sense} if either
there is a number $\eta\in(a,b)$ such that
$$\phantom{}_a^\eta\mathcal{K}(f) = \phantom{}_a^b\mathcal{K}(f), \quad
\textrm{or}\quad \lim_{x\to a^+} \phantom{}_a^x\mathcal{K}(f) =
\phantom{}_a^b\mathcal{K}(f).$$ Now, for $f\in\mathcal{C}\langle
a,b\rangle$ denote by $M$ the set of all points $x\in\langle
a,b\rangle$ in which $f$ is non-differentiable and put
$\mathfrak{m}=|M|$. Define the function
$$\mathcal{F}(x) :=
\frac{1}{x-a}\left(f'(x)-\phantom{}_a^x\mathcal{K}(f)\right),\quad
x\in(a,b\rangle\setminus M.$$ Then the assertion of Flett's
theorem is equivalent to $\mathcal{F}(\eta)=0$. Clearly, if
$\mathfrak{m}=0$, then by Riedel-Sahoo's theorem there exists
$\eta\in(a,b)$ such that
$$\mathcal{F}(\eta)=\frac{1}{2}\cdot\phantom{}_a^b\mathcal{K}(f').$$
So, what if $\mathfrak{m}>0$?

\begin{theorem}[Powers-Riedel-Sahoo, 2001]\label{thmFl}
Let $f\in\mathcal{C}\langle a,b\rangle$.
\begin{itemize}
\item[(i)] If $\mathfrak{m}\leq n$ for some non-negative integer
and $a\notin M$, then there exist $n+1$ points
$\eta_1,\dots,\eta_{n+1}\in(a,b)$ and $n+1$ positive numbers
$\alpha_1,\dots,\alpha_{n+1}$ with $\sum_{i=1}^{n+1}\alpha_i = 1$
such that
$$\sum_{i=1}^{n+1}\alpha_i\mathcal{F}(\eta_i) = \frac{1}{b-a}\left(\phantom{}_a^b\mathcal{K}(f)-f'(a)\right).$$
\item[(ii)] If $\mathfrak{m}$ is infinite and the graph of $f$
intersects its chord in the extended sense, then there exist
$\eta\in (a,b)$ and two positive numbers $\delta_1, \delta_2$ such
that
$$\textrm{either}\quad\mathcal{F}_1(\eta,h)\leq 0\leq \mathcal{F}_2(\eta,k), \quad \textrm{or}\quad
\mathcal{F}_2(\eta,k)\leq 0\leq \mathcal{F}_1(\eta,h),$$ holds for
$h\in(0,\delta_1\rangle$ and $k\in(0,\delta_2\rangle$, where
$$\mathcal{F}_1(\eta,h):=(\eta-a)\left(\phantom{}_{\eta-h}^\eta\mathcal{K}(f)-\phantom{}_a^\eta\mathcal{K}(f)\right),$$
$$\mathcal{F}_2(\eta,k):=(\eta-a)\left(\phantom{}_{\eta}^{\eta+k}\mathcal{K}(f)-\phantom{}_a^\eta\mathcal{K}(f)\right).$$
\end{itemize}
\end{theorem}

In item (i) we note that if $f'(a) =
\phantom{}_a^b\mathcal{K}(f)$, i.e., the second condition for the
graph of $f$ intersecting its chord in the extended sense holds,
then the convex combination of values of $\mathcal{F}$ at points
$\eta_i$, $i=1,\dots,n+1$, is simply zero. If, in item (ii), $f$
is differentiable at $\eta$, then $$\lim_{h\to 0^+}
\frac{\mathcal{F}_1(\eta,h)}{(\eta-a)^2} = \lim_{k\to 0^+}
\frac{\mathcal{F}_2(\eta,k)}{(\eta-a)^2} = \mathcal{F}(\eta).$$
The proof of item (i) can be found in~\cite{PRS01} and the proof
of (ii) is given in~\cite{DV}. Note that in the paper~\cite{PRS01}
authors extended the results of Theorem~\ref{thmFl} in the context
of topological vector spaces $X, Y$ for a class of Gateaux
differentiable functions $f:X\to Y$.

Flett's and Riedel-Sahoo's theorem give an opportunity to study
the behaviour of intermediate points from different points of
view. Recall that points $\eta$ (depending on the interval
$\langle a,b\rangle$) from Flett's, or Riedel-Sahoo's theorem are
called \textit{the Flett's}, or \textit{the Riedel-Sahoo's points
of function $f$ on the interval $\langle a,b\rangle$},
respectively.

The questions of stability of Flett's points was firstly
investigated in~\cite{DRS03}, but the main result therein was
shown to be incorrect. In paper~\cite{LXY} the correction was made
and the following results on Hyers-Ulam's stability of
Riedel-Sahoo's and Flett's points were proved.

\begin{theorem}[Lee-Xu-Ye, 2009]
Let $f\in\mathcal{D}\langle a, b\rangle$ and $\eta$ be a
Riedel-Sahoo's point of $f$ on $\langle a, b\rangle$. If $f$ is
twice differentiable at $\eta$ and
$$f''(\eta)(\eta-a)-2f'(\eta)+ 2\,\phantom{}_a^\eta\mathcal{K}(f)\neq
0,$$ then to any $\varepsilon>0$ and any neighborhood $N \subset
(a, b)$ of $\eta$, there exists a $\delta > 0$ such that for every
$g\in\mathcal{D}\langle a, b\rangle$ satisfying $|g(x)-
g(a)-(f(x)- f(a))| < \delta$ for $x\in N$ and $g'(b)-g'(a) =
f'(b)-f'(a)$, there exists a point $\xi \in N$ such that $\xi$ is
a Riedel-Sahoo's point of $g$ and $|\xi-\eta| < \varepsilon$.
\end{theorem}

As a corollary we get the Hyers-Ulam's stability of Flett's
points.

\begin{theorem}[Lee-Xu-Ye, 2009]
Let $f\in\mathcal{D}\langle a, b\rangle$ with $f'(a)=f'(b)$ and
$\eta$ be a Flett's point of $f$ on $\langle a, b\rangle$. If $f$
is twice differentiable at $\eta$ and
$$f''(\eta)(\eta-a)-2f'(\eta)+ 2\,\phantom{}_a^\eta\mathcal{K}(f)\neq
0,$$ then to any $\varepsilon>0$ and any neighborhood $N \subset
(a, b)$ of $\eta$, there exists a $\delta > 0$ such that for every
$g\in\mathcal{D}\langle a, b\rangle$ satisfying $g(a)=f(a)$ and
$|g(x)-f(x)| < \delta$ for $x\in N$, there exists a point $\xi \in
N$ such that $\xi$ is a Flett's point of $g$ and $|\xi-\eta| <
\varepsilon$.
\end{theorem}

Another interesting question is the limit behaviour of
Riedel-Sahoo's points (Flett's points are not interesting because
of the condition $f'(a)=f'(b)$). We demonstrate the main idea on
the following easy example: let $f(t)=t^3$ for $t\in\langle
0,x\rangle$ with $x>0$. By Riedel-Sahoo's theorem for each $x>0$
there exists a point $\eta_x\in(0,x)$ such that
$$3\eta_x^2 =
\frac{\eta_x^3}{\eta_x}+\frac{3x^2}{x}\cdot\frac{\eta_x}{2}\,\,\,
\Leftrightarrow \,\,\, 4\eta_x^2 = 3x\eta_x\,\,\, \Leftrightarrow
\,\,\, \eta_x = \frac{3}{4}x.$$ Thus, we have obtained a
dependence of Riedel-Sahoo's points on $x$. If we shorten the
considered interval, we get
$$\lim_{x\to 0^+} \frac{\eta_x-0}{x-0} = \lim_{x\to 0^+}
\frac{\frac{3}{4}x}{x} = \frac{3}{4}.$$ So, how do Flett's points
behave for the widest class of function? In paper~\cite{PRS98}
authors proved the following result.

\begin{theorem}[Powers-Riedel-Sahoo, 1998]
Let $f\in\mathcal{D}\langle a,a+x\rangle$ be such that
$$f(t) = p(t)+(t-a)^\alpha g(t), \quad
\alpha\in(1,2)\cup(2,+\infty),
$$ where $p$ is a polynomial at most of second order, $g'$ is
bounded on the interval $(a,a+x\rangle$ and
$g(a)=\lim\limits_{x\to 0^+} g(a+x)\neq 0$. Then
$$\lim_{x\to 0^+} \frac{\eta_x-a}{x} =
\left(\frac{\alpha}{2(\alpha-1)}\right)^{\frac{1}{\alpha-2}},$$
where $\eta_x$ are the corresponding Riedel-Sahoo's points of $f$
on $\langle a,a+x\rangle$.
\end{theorem}

\begin{problem}
Enlarge the Power-Riedel-Sahoo's family of functions for which it
is possible to state the exact formula for limit properties of
corresponding intermediate points.
\end{problem}

\section{Further sufficient conditions for validity
of~(\ref{flettMV})}\label{sectionsufficient}

In this section we review some other conditions yielding validity
of equality~(\ref{flettMV}).

\subsection{Trahan's inequalities}

Probably the first study about Flett's result and its
generalization is dated to the year 1966 by~\textsc{Donald H.
Trahan}~\cite{trahan}. He provides a different condition for the
assertion of Flett's theorem under some inequality using a
comparison of slopes of secant line passing through the end points
and tangents at the end points.

\begin{theorem}[Trahan, 1966]\label{trahan}
Let $f\in\mathcal{D}\langle a,b\rangle$ and
\begin{equation}\label{predpoklad}
\left(f'(b)-\phantom{}_a^b\mathcal{K}(f)\right)\cdot
\left(f'(a)-\phantom{}_a^b\mathcal{K}(f)\right)\geq 0.
\end{equation}
Then there exists $\eta\in(a,b\rangle$ such that~(\ref{flettMV})
holds.
\end{theorem}

Donald Trahan in his proof again considers the function $g$ given
by~(\ref{flettfunkcia}). Then $g\in\mathcal{C}\langle
a,b\rangle\cap\mathcal{D}(a,b\rangle$ and
$$g'(x)=\frac{1}{x-a}\left(f'(x)-\phantom{}_a^x\mathcal{K}(f)\right), \quad x\in(a,b\rangle.$$
Since
$$[g(b)-g(a)]\,g'(b)=\frac{-1}{b-a}\left(f'(b)-\phantom{}_a^b\mathcal{K}(f)\right)\cdot\left(f'(a)-\phantom{}_a^b\mathcal{K}(f)\right),$$
then by~(\ref{predpoklad}) we get $[g(b)-g(a)]\,g'(b)\leq 0$. Now
Trahan concludes that $g'(\eta)=0$ for some $\eta\in(a,b\rangle$,
which is equivalent to~(\ref{flettMV}).

The only step here is to prove \textsf{Trahan's lemma}, i.e., the
assertion of Rolle's theorem under the conditions
$g\in\mathcal{C}\langle a,b\rangle\cap\mathcal{D}(a,b\rangle$ and
$[g(b)-g(a)]\,g'(b)\leq 0$. Easily, if $g(a)=g(b)$, then Rolle's
theorem gives the desired result. If $g'(b)=0$, putting $\eta=b$
we have $g'(\eta)=0$. So, let us assume that
$[g(b)-g(a)]\,g'(b)<0$. This means that either $g'(b)<0$ and
$g(b)>g(a)$, or $g'(b)>0$ and $g(b)<g(a)$. In the first case,
since $g\in\mathcal{C}\langle a,b\rangle$, $g(b)>g(a)$ and $g$ is
strictly decreasing in $b$, then $g$ has its maximum at
$\eta\in(a,b)$ and by Fermat's theorem we get $g'(\eta)=0$.
Similarly, in the second case $g$ has minimum at the same point
$\eta\in(a,b)$, thus $g'(\eta)=0$.

\begin{remark}\rm\label{remarkTrahanFlett}
Obviously, the class of Trahan's functions, i.e., differentaible
functions on $\langle a,b\rangle$ satisfying Trahan's
condition~(\ref{predpoklad}), is wider than the class of Flett's
functions $f\in\mathcal{D}\langle a,b\rangle$ satisfying Flett's
condition $f'(a)=f'(b)$. Indeed, for $f'(a)=f'(b)$ Trahan's
condition~(\ref{predpoklad}) is trivially fulfilled. On the other
hand the function $y=x^3$ for $x\in\langle -\frac{1}{2},1\rangle$
does not satisfy Flett's condition, and it is easy to verify that
it satisfies Trahan's one.
\end{remark}

\begin{figure}
\begin{center}
\includegraphics[scale=0.8]{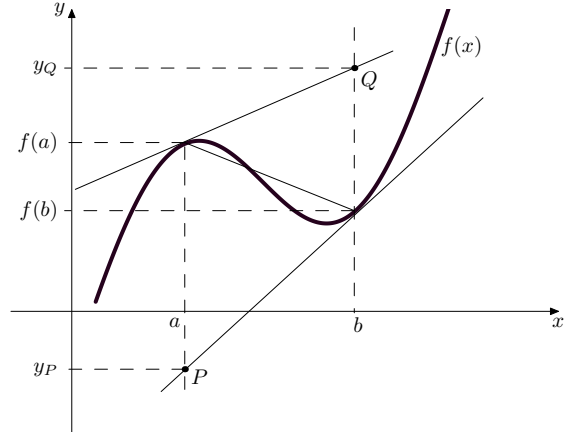}
\caption{Geometrical interpretation of Trahan's
condition}\label{obrTrahan}
\end{center}
\end{figure}

\paragraph{Geometrical meaning of Trahan's condition} Clearly, Trahan's
inequality~(\ref{predpoklad}) holds if and only if
$$\left[f'(b)\geq\phantom{}_a^b\mathcal{K}(f) \wedge f'(a)\geq\phantom{}_a^b\mathcal{K}(f)\right]
\vee \left[f'(b)\leq\phantom{}_a^b\mathcal{K}(f) \wedge
f'(a)\leq\phantom{}_a^b\mathcal{K}(f)\right].$$ Since
$\phantom{}_a^b\mathcal{K}(f)$ gives the slope of the secant line
between $[a,f(a)]$ and $[b,f(b)]$, Trahan's condition requires
either both slopes of tangents at the end points are greater or
equal, or both are smaller or equal to
$\phantom{}_a^b\mathcal{K}(f)$. We consider two cases:
\begin{itemize}
\item[(i)] if $f'(b)=\phantom{}_a^b\mathcal{K}(f)$, then the
tangent at $b$ is parallel to the secant, and the tangent at $a$
may be arbitrary (parallel to the secant, lying above or under the
graph of secant on $(a,b)$), analogously for
$f'(a)=\phantom{}_a^b\mathcal{K}(f)$;

\item[(ii)] if $f'(b)\neq\phantom{}_a^b\mathcal{K}(f)$ and
$f'(a)\neq\phantom{}_a^b\mathcal{K}(f)$, then one of the tangents
at the end points has to lie above and the second one under the
graph of secant line on $(a,b)$, or vice versa, see
Fig.~\ref{obrTrahan}. More precisely, let tangent at $a$ intersect
the line $x=b$ at the point $Q=[b,y_Q]$ and tangent at $b$
intersect the line $x=a$ at the point $P=[a,y_P]$. Then either
$y_Q>f(b)$ and $y_P<f(a)$, or $y_Q<f(b)$ and $y_P>f(a)$. For
parallel tangents at the end points, i.e., for $f'(a)=f'(b)$, this
geometrical interpretation provides a new insight which leads to
the already mentioned paper~\cite{DV}.
\end{itemize}

Moreover, Trahan in his paper~\cite{trahan} provides other
generalization of Flett's theorem. Namely, he proves certain
,,Cauchy form`` of his result for two functions which will be a
source of our results later in Section~\ref{secvlastne}.

\begin{theorem}[Trahan, 1966]\label{trahveta*}
Let $f, g\in\mathcal{D}\langle a,b\rangle$, $g'(x)\neq 0$ for each
$x\in\langle a,b\rangle$ and
$$\left(\frac{f'(a)}{g'(a)}-\phantom{}_a^b\mathcal{K}(f,
g)\right)\Bigl(\phantom{}_a^b\mathcal{K}(g)\,
f'(b)-\phantom{}_a^b\mathcal{K}(f)\,g'(b)\Bigr)\geq 0.$$ Then
there exists $\eta\in(a,b\rangle$ such that
%{\setlength\arraycolsep{2pt}\begin{eqnarray}\label{Trvztah}
$\displaystyle
\frac{f'(\eta)}{g'(\eta)}=\phantom{}_a^\eta\mathcal{K}(f, g).
$%\end{eqnarray}}
\end{theorem}

Its proof is based on application of Trahan's
lemma~\cite[Lemma~1]{trahan} for function
$$h(x) =
\begin{cases}
\phantom{}_a^x\mathcal{K}(f,g), & x\in(a,b\rangle \\
\frac{f'(a)}{g'(a)}, & x=a.
\end{cases}
$$

%A different sufficient condition for the validity
%of~(\ref{Trvztah}) is as follows.
%\begin{proposition}
%Let $f, g\in\mathcal{D}\langle a,b\rangle$ such that $g'(x)\neq 0$
%for each $x\in\langle a,b\rangle$. If
%$\phantom{}_a^b\mathcal{K}(g)\,g'(b)>0$ and
%$\frac{f'(a)}{g'(a)}=\frac{f'(b)}{g'(b)}$, then there exists
%$\eta\in(a,b\rangle$ such that (\ref{Trvztah}) holds.
%\end{proposition}

\subsection{Tong's discrete and integral means}

Another sufficient condition for validity of~(\ref{flettMV}) was
provided by \textsc{JingCheong Tong} in the beginning of 21st
century in his paper~\cite{tong}. Tong does not require
differentiability of function $f$ at the end points of the
interval $\langle a,b\rangle$, but he uses certain means of that
function. Indeed, for a function $f: M\to \mathbb{R}$ and two
distinct points $a,b\in M$ denote by {\setlength\arraycolsep{2pt}
\begin{eqnarray*}\label{priemer}
A_f(a,b) = \frac{f(a)+f(b)}{2}\quad\textrm{and}\quad I_f(a,b) =
\frac{1}{b-a}\int_{a}^{b}f(t)\,\mathrm{d}t
\end{eqnarray*}}the \textit{arithmetic} (discrete) and \textit{integral} (continuous) \textit{mean
of $f$ on the interval} $\langle a,b\rangle$, respectively.

\begin{theorem}[Tong, 2004]\label{tong}
Let $f\in\mathcal{C}\langle a,b\rangle\cap\mathcal{D}(a,b)$. If
$A_f(a,b)=I_f(a,b)$, then there exists $\eta\in(a,b)$ such that
(\ref{flettMV}) holds.
\end{theorem}

In his proof Tong defines the auxiliary function
$$h(x)=\begin{cases}
(x-a)[A_f(a,x)-I_f(a,x)], & x\in (a,b\rangle \\ 0, & x=a.
\end{cases}$$ Easily, $h\in\mathcal{C}\langle
a,b\rangle\cap\mathcal{D}(a,b)$ and $h(a)=0=h(b)$. Then Rolle's
theorem for $h$ on $\langle a,b\rangle$ finishes the proof.

\paragraph{Geometrical meaning of Tong's condition}
The condition $A_f(a,b)=I_f(a,b)$ is not so evident geometrically
in comparison with the Flett's condition $f'(a)=f'(b)$. In some
sense we can demonstrate it as "the area under the graph of $f$ on
$\langle a,b\rangle$ is exactly the volume of a rectangle with
sides $b-a$ and $\frac{f(a)+f(b)}{2}$``.%, see Fig.???.

Let us analyze Tong's condition $A_f(a,b)=I_f(a,b)$ for
$f\in\mathcal{C}\langle a,b\rangle\cap\mathcal{D}(a,b)$ in detail.
It is important to note that this equality does not hold in
general for each $f\in\mathcal{C}\langle
a,b\rangle\cap\mathcal{D}(a,b)$. Indeed, for $f(x)=x^{2}$ on
$\langle 0,1\rangle$ we have
$$A_f(0,1)=\frac{0^{2}+1^{2}}{2}=\frac{1}{2}, \qquad I_f(0,1)=\frac{1}{1-0}\int_{0}^{1} x^2\,\mathrm{d}x=\frac{1}{3}.$$
A natural question is \textit{how large is the class of such
functions?} For $f\in\mathcal{C}(M)\cap\mathcal{D}(M)$ denote by
$F$ a primitive function to $f$ on an interval $M$ and let $a,b$
be interior points of $M$. Then the condition $A_f(x,b)=I_f(x,b)$,
$x\in M$, is equivalent to the condition
\begin{equation}\label{fF}
\frac{f(x)+f(b)}{2} = \phantom{}_a^x\mathcal{K}(F), \quad x\neq
a.\end{equation} Since $f\in\mathcal{D}(M)$, then
$F\in\mathcal{D}^2(M)$ and $f'(t)=F''(t)$ for each $t\in M$.
Differentiating the equality~(\ref{fF}) with respect to $x$ we get
$$\frac{f'(x)}{2} = \frac{F'(x)-\phantom{}_a^x\mathcal{K}(F)}{x-a},$$ which
is equivalent to the equation
$$F''(x)(x-a)^2 = 2(F'(x)(x-a)+F(a)-F(x)).$$ Solving this
differential equation on intervals $(-\infty, b)\cap M$ and
$(b,+\infty)\cap M$, and using the second differentiability of $F$
at $b$ we have
$$F(x)=\frac{\alpha}{2}x^2+\beta x+\gamma, \quad x\in M,
\,\,\alpha, \beta, \gamma \in\mathbb{R}, \alpha\neq 0,$$ and
therefore
$$f(x)=\alpha x+\beta, \quad x\in M.$$ So, the class of functions fulfilling Tong's condition $A_f(a,b)=I_f(a,b)$
for each interval $\langle a,b\rangle$ is quite small (affine
functions, in fact). Of course, if we do not require the condition
"on each interval $\langle a,b\rangle$", then we may use, e.g.,
the function $y=\arcsin x$ on the interval $\langle -1,1\rangle$
which does not satisfy neither Flett's nor Trahan's condition
(because it is not differentiable at the end points).

\begin{figure}
\begin{center}
\includegraphics[scale=0.8]{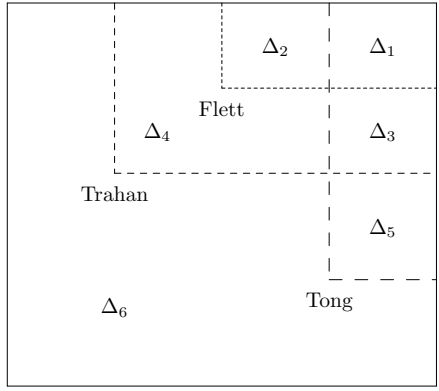}
\caption{The relations among Flett's, Trahan's and Tong's families
of functions}\label{obrdiagramy}
\end{center}
\end{figure}

\begin{remark}\rm\label{rempriklady}
The relations among the classes of Flett's, Trahan's and Tong's
functions are visualized in Fig.~\ref{obrdiagramy}, where each
class is displayed as a rectangle with the corresponding name
below the left corner. Moreover, $\Delta_i$, $i=1,\dots, 6$, are
the classes of functions of possible relationships of Flett's,
Trahan's and Tong's classes of functions. For instance, $\Delta_6$
denotes a class of (not necessarily differentiable or continuous)
functions on $\langle a,b\rangle$ for which none of the three
conditions is fulfilled, but the assertion of Flett's theorem
still holds. Immediately, $\Delta_1$ is non-empty, because it
contains all affine functions on $\langle a,b\rangle$. Thus,
\begin{itemize}
\item[(i)] Flett's and Trahan's conditions were compared in
Remark~\ref{remarkTrahanFlett} yielding that Trahan's class of
functions is wider than Flett's one, i.e.,
$\Delta_1\cup\Delta_2\subset\Delta_3\cup\Delta_4$; \item[(ii)]
Tong's condition and Flett's condition are independent each other,
because for
$$f(x)=\sin x, \quad x\in\left\langle
\frac{\pi}{2},\frac{5}{2}\pi\right\rangle,$$ we have
$f'(\frac{\pi}{2})=f'(\frac{5}{2}\pi)=0$, but
$1=A_f(\frac{\pi}{2},\frac{5}{2}\pi)\neq
I_f(\frac{\pi}{2},\frac{5}{2}\pi)=0$; on the other hand
$f(x)=\arcsin x$ for $x\in\langle-1,1\rangle$ fulfills Tong's
condition, but does not satisfy Trahan's one; also for $f(x)=x^3$,
$x\in\langle -1,1\rangle$, we have $A_f(-1,1)=I_f(-1,1)$ and
$f'(-1)=f'(1)$ which yields that the classes of functions
$\Delta_1$, $\Delta_2$ and $\Delta_3$ are non-empty; \item[(iii)]
similarly, Trahan's condition and Tong's condition are
independent, e.g., the function $f(x)=x^3$ on the interval
$\langle -\frac{1}{2},1\rangle$ satisfies Trahan's condition, but
$A_f(-\frac{1}{2},1)\neq I_f(-\frac{1}{2},1)$, so $\Delta_4$ is
non-empty, and for $f(x)=\arcsin x$, $x\in\langle -1,1\rangle$ we
have $f\in\Delta_5$; \item[(iv)] for the function $\textrm{sgn}$
on $\langle -2,1\rangle$ none of the three conditions is
fulfilled, but the assertion~(\ref{flettMV}) still holds, i.e.,
$\Delta_6$ is non-empty.
\end{itemize}
\end{remark}

\begin{question}
Is each class $\Delta_i$, $i\in\{3,5,6\}$, non-empty when
considering the stronger condition $f\in\mathcal{D}\langle
a,b\rangle$ in Tong's assumption?
\end{question}

Removing the condition $A_f(a,b)=I_f(a,b)$ Tong obtained the
following result which no more corresponds to the result of
Riedel-Sahoo's theorem.

\begin{theorem}[Tong, 2004]\label{tong2}
Let $f\in\mathcal{C}\langle a,b\rangle\cap\mathcal{D}(a,b)$. Then
there exists $\eta\in(a,b)$ such that
$$f'(\eta)=\phantom{}_a^\eta\mathcal{K}(f)+\frac{6[A_f(a,b)-I_f(a,b)]}{(b-a)^{2}}(\eta-a).$$
\end{theorem}

Tong's proof uses the auxiliary function
$$H(x)=f(x)-\frac{6[A_f(a,b)-I_f(a,b)]}{(b-a)^{2}}(x-a)(x-b), \quad x\in \langle a,b\rangle.$$
It is easy to verify that $H\in\mathcal{C}\langle
a,b\rangle\cap\mathcal{D}(a,b)$, $H(a)=f(a)$ and $H(b)=f(b)$.
Thus, $A_H(a,b)=A_f(a,b)$ and $I_H(a,b) = A_H(a,b)$. Then by
Theorem~\ref{tong} there exists $\eta\in(a,b)$ such that
$H'(\eta)=\phantom{}_a^\eta\mathcal{K}(H)$ which is equivalent to
the assertion of theorem.

\begin{question}
Analogously to Riedel-Sahoo's Theorem~\ref{thmRS} we may ask the
following: What is the limit behaviour of Tong's points $\eta$ of
a function $f$ on the interval $\langle a,b\rangle$?
\end{question}

\subsection{Two new extensions of Flett's
theorem}\label{secvlastne}

In this section we present other sufficient conditions for
validity of~(\ref{flettMV}) and its extension. As far as we know
they are not included in any literature we were able to find. The
basic idea is a mixture of Trahan's results with (although not
explicitly stated) Diaz-V\'yborn\'y's concept of intersecting the
graphs of two functions~\cite{DV}. We will also present nice
geometrical interpretations of these results. A particular case of
our second result is discussed in the end of this section.

\begin{lemma}\label{lema1}
If $f, g\in\mathcal{D}\langle a,b\rangle$, $g(b)\neq g(a)$ and
\begin{equation}\label{novapod0}
\left[f'(a)-\phantom{}_a^b\mathcal{K}(f,g)
\,g'(a)\right]\cdot\left[f'(b)-\phantom{}_a^b\mathcal{K}(f,g)\,g'(b)\right]\geq
0,\end{equation} then there exists $\xi\in(a,b)$ such that
\begin{equation}\label{novapod}
f(\xi)-f(a)=\phantom{}_a^b\mathcal{K}(f,g)(g(\xi)-g(a)).
\end{equation}
\end{lemma}

\proof Let us consider the function
$$\varphi(x)=f(x)-f(a)-\phantom{}_a^b\mathcal{K}(f,g)\cdot(g(x)-g(a)), \quad x\in \langle a,b\rangle.$$
Then $$\varphi'(a) = f'(a)-\phantom{}_a^b\mathcal{K}(f,g)\cdot
g'(a) \quad \textrm{and} \quad \varphi'(b) =
f'(b)-\phantom{}_a^b\mathcal{K}(f,g)\cdot g'(b).$$ If
$\varphi'(a)\geq 0$, then according to assumption we get
$\varphi'(b)\geq 0$. So, there exist points $\alpha,\beta\in(a,b)$
such that $\varphi(\alpha)>0$ and $\varphi(\beta)<0$. Thus,
$\varphi(\alpha)\varphi(\beta)<0$ and by Bolzano's theorem there
exists a point $\xi\in(\alpha,\beta)$ such that $\varphi(\xi)=0$.
The case $\varphi'(a)\leq 0$ and $\varphi'(b)\leq 0$ is analogous.
\qed

\begin{theorem}\label{thmnew1}
If $f, g\in\mathcal{D}\langle a,b\rangle$, $g(b)\neq g(a)$ and the
condition~(\ref{novapod0}) holds, then there exists $\eta\in(a,b)$
such that
\begin{equation}\label{novyvys}
f'(\eta)-\phantom{}_a^\eta\mathcal{K}(f) =
\phantom{}_a^b\mathcal{K}(f,g)\left[g'(\eta)-\phantom{}_a^\eta\mathcal{K}(g)\right].
\end{equation}
\end{theorem}

\proof Let us take the auxiliary function
$$
F(x) = \begin{cases}
\phantom{}_a^x\mathcal{K}(f) - \phantom{}_a^b\mathcal{K}(f,g)\cdot\phantom{}_a^x\mathcal{K}(g), & x\in (a,b\rangle,\\
f'(a) - \phantom{}_a^b\mathcal{K}(f,g)\cdot g'(a), & x=a.
\end{cases}
$$
Observe that $F(x)=\phantom{}_a^x\mathcal{K}(\varphi)$ for
$x\in(a,b\rangle$, where $\varphi$ is the auxiliary function from
the proof of Lemma~\ref{lema1}. Thus, by Lemma~\ref{lema1} there
exists a point $\xi\in(a,b)$ such that $F(\xi)=0=F(b)$. Since
$F\in\mathcal{C}\langle \xi,b\rangle\cap\mathcal{D}(\xi,b)$, then
by Rolle's theorem there exists $\eta\in(\xi,b)\subset(a,b)$ such
that $F'(\eta)=0$ which is equivalent to the desired result. \qed

\begin{figure}
\begin{center}
\includegraphics[scale=0.8]{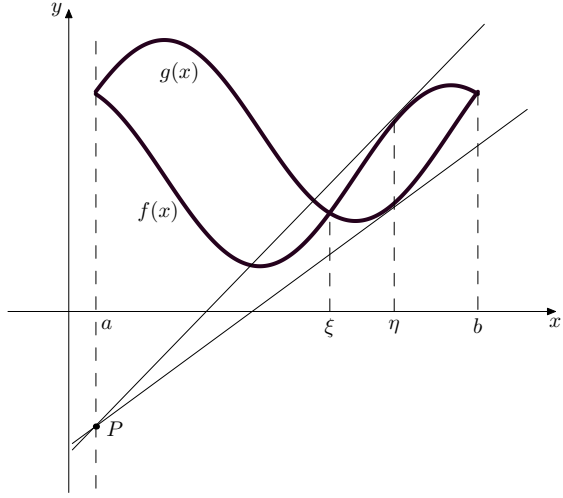}
\caption{Geometrical interpretation of Theorem~\ref{thmnew1}}
\label{obrFlett2}
\end{center}
\end{figure}

In what follows we denote by $$T_n(f,x_0)(x) :=
f(x_0)+\frac{f'(x_0)}{1!}(x-x_0)+\dots +
\frac{f^{(n)}(x_0)}{n!}(x-x_0)^n$$ the $n$-th Taylor's polynomial
of a function $f$ at a point $x_0$. Rewriting the assertion of
Theorem~\ref{thmnew1} in terms of Taylor's polynomial yields
$$f(a)-T_1(f,\eta)(a) =
\phantom{}_a^b\mathcal{K}(f,g)\cdot(g(a)-T_1(g,\eta)(a)).$$

\paragraph{Geometrical meaning of Theorem~\ref{thmnew1}}
Realize that $T_1(f,x_0)(x)$ is the tangent to the graph of $f$ at
the point $x_0$, i.e.,
$$T_1(f,x_0)(x)=f(x_0)+f'(x_0)(x-x_0).$$ Then the equation~(\ref{novyvys})
may be equivalently rewritten as follows
\begin{equation}\label{novyvys1}
(\exists \eta\in(\xi,b)\subset (a,b))\,\,f(a)-T_1(f,\eta)(a) =
\phantom{}_a^b\mathcal{K}(f,g)\cdot(g(a)-T_1(g,\eta)(a)).
\end{equation}
Since $$f(b)-f(a) = g(b)-g(a) \Leftrightarrow
f(b)-g(b)=f(a)-g(a),$$ then the equation~(\ref{novapod}) has the
form
\begin{equation}\label{novapod1}
(\exists \xi\in (a,b))\,\,f(\xi)-g(\xi) = f(a)-g(a) = f(b)-g(b),
\end{equation} and~(\ref{novyvys1}) may be rewritten into
\begin{equation}\label{novyvys2}
(\exists \eta\in(\xi,b))\,\,T_1(f,\eta)(a)-T_1(g,\eta)(a) =
f(a)-g(a) = f(b)-g(b).
\end{equation}
Thus, considering a function $g$ such that $g(a)=f(a)$ and
$g(b)=f(b)$ the equation~(\ref{novapod1}) yields
$$(\exists \xi\in (a,b))\,\, f(\xi) = g(\xi)$$ and
the equation~(\ref{novyvys2}) has the form
$$(\exists \eta\in (\xi,b))\,\, T_1(f,\eta)(a) = T_1(g,\eta)(a).$$
Geometrically it means that tangents at the points
$[\eta,f(\eta)]$ and $[\eta,g(\eta)]$ pass though the common point
$P$ on the line $x=a$, see Fig.~\ref{obrFlett2}.

\begin{figure}
\begin{center}
\includegraphics[scale=0.8]{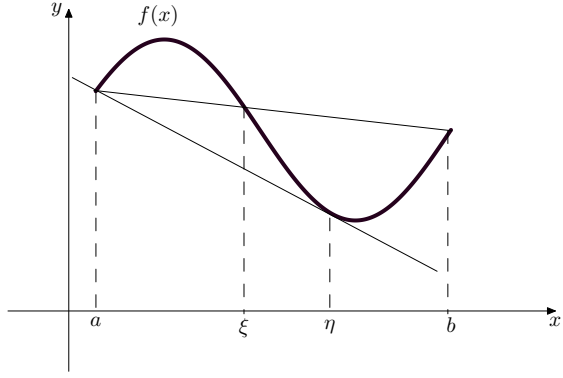}
\caption{Geometrical interpretation of Theorem~\ref{thmnew1} for
secant} \label{obrFlett3}
\end{center}
\end{figure}

\begin{remark}\rm
Observe that in the special case of secant joining the end points,
i.e., the function
$$g(x)= f(a)+\phantom{}_a^b\mathcal{K}(f)(x-a),$$ where $f$ is a
function fulfilling assumptions of Theorem~\ref{thmnew1}, we get
the original Trahan's result of Theorem~\ref{trahan} (in fact, a
generalization of Flett's theorem) with the explicit geometrical
interpretation on Fig.~\ref{obrFlett3}.
\end{remark}

\begin{lemma}\label{lema2}
Let $f, g\in\mathcal{D}\langle a,b\rangle$ and $f,g$ be twice
differentiable at the point $a$. If $g(a)\neq g(b)$ and
\begin{equation}\label{novapod00}
\left[f'(a)-\phantom{}_a^b\mathcal{K}(f,g)\cdot
g'(a)\right]\left[f''(a)-\phantom{}_a^b\mathcal{K}(f,g)\cdot
g''(a)\right]> 0,
\end{equation} then there exists $\xi\in(a,b)$
such that
$$f'(a)-\phantom{}_a^\xi\mathcal{K}(f)=\phantom{}_a^b\mathcal{K}(f,g)\left[g'(a)-\phantom{}_a^\xi\mathcal{K}(g)\right].$$
\end{lemma}

\proof Consider the function
$$
F(x) = \begin{cases} 0, & x=a, \\
f'(a) - \phantom{}_a^x\mathcal{K}(f)+\phantom{}_a^b\mathcal{K}(f,g)\left[\phantom{}_a^x\mathcal{K}(g)-g'(a)\right], & x\in (a,b),\\
f'(a) - \phantom{}_a^b\mathcal{K}(f,g)\cdot g'(a), & x=b.
\end{cases}$$ Then $F\in\mathcal{C}\langle a,b\rangle\cap
\mathcal{D}\langle a,b)$ with {\setlength\arraycolsep{2pt}
\begin{eqnarray}\label{derivacia}F'(a) & = & %\lim_{x\to a^+}\phantom{}_a^x\mathcal{K}(F) =
\lim_{x\to a^+}\frac{f'(a)(x-a)-(f(x)-f(a))+\phantom{}_a^b\mathcal{K}(f,g)\left[g(x)-g(a)-g'(a)(x-a)\right]}{(x-a)^2}\nonumber\\
& = &
\lim_{x\to a^+}\frac{f'(a)-f'(x)+\phantom{}_a^b\mathcal{K}(f,g)\cdot [g'(x)-g'(a)]}{2(x-a)}\nonumber\\
& = & -\frac{1}{2}\lim_{x\to
a^+}\left(\frac{f'(x)-f'(a)}{x-a}-\phantom{}_a^b\mathcal{K}(f,g)\cdot
\frac{g'(x)-g'(a)}{x-a}\right)\nonumber \\ & = &
-\frac{1}{2}\left[f''(a)-\phantom{}_a^b\mathcal{K}(f,g)\cdot
g''(a)\right],\nonumber
\end{eqnarray}}where L'Hospital rule has been used. Suppose that
$$F(b)=\left[f'(a)-\phantom{}_a^b\mathcal{K}(f,g)\cdot
g'(a)\right]>0,$$ analogous arguments apply if $F(b)<0$. Then
$$\left[f''(a)-\phantom{}_a^b\mathcal{K}(f,g)\cdot
g''(a)\right]>0$$ by assumption of theorem which implies
$F'(a)<0$. Since $F(a)=0$, then there exists $\alpha\in(a,b)$ such
that $F(\alpha)<0$. According to Bolzano's theorem there exists a
point $\xi\in(\alpha,b)$ such that $F(\xi)=0$, which completes the
proof. \qed

\begin{theorem}\label{veta2}
Let $f, g\in\mathcal{D}\langle a,b\rangle$ and $f,g$ be twice
differentiable at the point $a$. If $g(a)\neq g(b)$ and the
inequality~(\ref{novapod00}) holds, then there exists
$\eta\in(a,b)$ such that~(\ref{novyvys}) holds.
\end{theorem}

\proof Consider the function $F$ as in the proof of
Lemma~\ref{lema2}. Then by Lemma~\ref{lema2} there exists
$\xi\in(a,b)$ such that $F(\xi)=0=F(a)$ and by Rolle's theorem
there exists $\eta\in(a,\xi)$ such that $F'(\eta)=0$. \qed

\begin{figure}
\begin{center}
\includegraphics[scale=0.8]{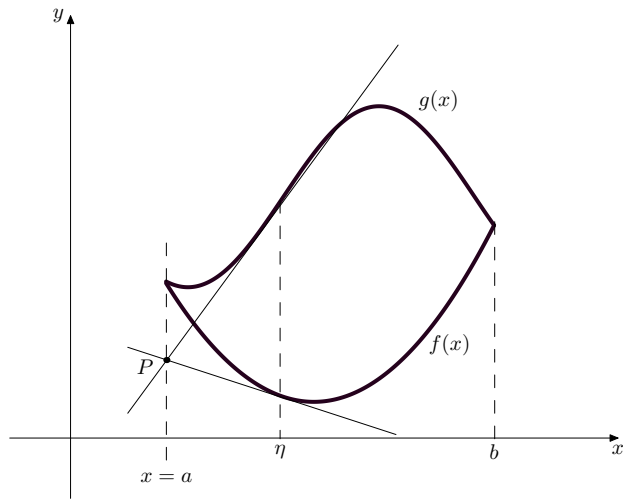}
\caption{Geometrical interpretation of Theorem~\ref{veta2}}
\label{obrFlett4}
\end{center}
\end{figure}

\paragraph{Geometrical meaning of Theorem~\ref{veta2}} Using the Taylor's polynomial we can rewrite the assertion of
Lemma~\ref{lema2} and Theorem~\ref{veta2} as follows
\begin{equation}\label{novyvys3}
(\exists \xi\in (a,b))\,\, f(\xi)-T_1(f,a)(\xi) =
\phantom{}_a^b\mathcal{K}(f,g)\cdot (g(\xi)-T_1(g,a)(\xi))
\end{equation}
and
\begin{equation}\label{novyvys4}
(\exists \eta\in (a,\xi))\,\, f(a)-T_1(f,a)(\eta) =
\phantom{}_a^b\mathcal{K}(f,g)\cdot (g(a)-T_1(g,a)(\eta)),
\end{equation} respectively.
Since $$f(b)-f(a)=g(b)-g(a) \,\,\Leftrightarrow\,\,
f(b)-g(b)=f(a)-g(a),$$ then~(\ref{novyvys3}) may be rewritten as
$$(\exists \xi\in (a,b))\,\, f(\xi)-g(\xi)=T_1(f,a)(\xi)-
T_1(g,a)(\xi).
$$ Similarly~(\ref{novyvys4}) may be rewritten as follows
$$(\exists \eta\in (a,\xi))\,\, f(a)-g(a) = T_1(f,\eta)(a)-T_1(g,\eta)(a)).
$$
If $f$ and $g$ have the same values at the end points, the last
equation reduces to
$$(\exists \eta\in (a,b))\,\, T_1(f,\eta)(a) = T_1(g,\eta)(a).$$ Geometrically it
means that tangents at points $[\eta,f(\eta)]$ and
$[\eta,g(\eta)]$ pass through the common point $P$ on the line
$x=a$, see Fig.~\ref{obrFlett4}.

\begin{figure}
\begin{center}
\includegraphics[scale=0.8]{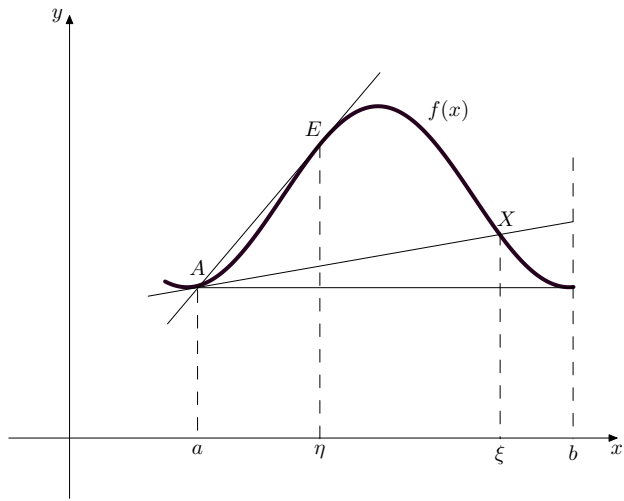}
\caption{Geometrical interpretation of Theorem~\ref{veta2} for
secant} \label{obrFlett5}
\end{center}
\end{figure}

\begin{remark}\rm
Again, Theorem~\ref{veta2} for the secant
$$g(x)= f(a)+\phantom{}_a^b\mathcal{K}(f)(x-a)$$ guarantees the existence of a point $\eta\in(a,b)$
such that $T_1(f,\eta)(a)=f(a)$, i.e., tangent at $[\eta,f(\eta)]$
passes through the point $A=[a,f(a)]$ which is exactly the
geometrical interpretation of Flett's theorem in
Fig.~\ref{obrFlett1}. The assumption~(\ref{novapod00}) reduces in
the secant case to the inequality
\begin{equation}\label{malesevic}
\left[f'(a)-\phantom{}_a^b\mathcal{K}(f)\right]f''(a)>0,
\end{equation}
i.e.,
$$\left[f'(a)>\phantom{}_a^b\mathcal{K}(f) \wedge f''(a)>0\right] \vee \left[f'(a)<\phantom{}_a^b\mathcal{K}(f) \wedge f''(a)<0\right].$$
Considering the first case yields {\setlength\arraycolsep{2pt}
\begin{eqnarray*}
f'(a)>\phantom{}_a^b\mathcal{K}(f) \wedge f''(a)>0 \,\, &
\Leftrightarrow & \,\, f(b) < f(a)+f'(a)(b-a) \wedge f''(a)>0 \\ &
\Leftrightarrow & \,\, f(b)<T_1(f,a)(b) \wedge f''(a)>0.
\end{eqnarray*}}This means that there exists a point
$X=[\xi,f(\xi)]$ such that the line $AX$ is a tangent to the graph
of $f$ at $A=[a,f(a)]$. Then from the assertion of
Theorem~\ref{veta2} we have the existence of a point
$E=[\eta,f(\eta)]$, where $\eta\in(a,\xi)$, such that the tangent
to the graph of $f$ at $E$ passes through the point $A=[a,f(a)]$,
see Fig.~\ref{obrFlett5}. Similarly for the second case.
\end{remark}

\begin{remark}\rm
Observe that if $f'(a)=\phantom{}_a^b\mathcal{K}(f)$, then the
condition~(\ref{malesevic}) is not fulfilled, but the
assertion~(\ref{flettMV}) of Flett's theorem still holds by
Trahan's condition. On the other hand, if
$f'(a)\neq\phantom{}_a^b\mathcal{K}(f)$ and $f''(a)=0$, then the
assertion~(\ref{flettMV}) does not need to hold, e.g., for
$f(x)=\sin x$ on the interval $\langle 0,\pi\rangle$ we have
$(f'(0)-\phantom{}_0^\pi\mathcal{K}(f))\cdot f''(0) = 1\cdot 0 =
0$, but there is no such a point $\eta\in(0,\pi)$ which is a
solution of the equation $\eta \cos\eta = \sin\eta$.
\end{remark}

We have to point out that the inequality~(\ref{malesevic}) was
observed as a sufficient condition for validity of~(\ref{flettMV})
in~\cite{Malesevic}, but starting from a different point,
therefore our general result of Theorem~\ref{veta2} seems to be
new. Indeed, Male\v{s}evi\`{c} in~\cite{Malesevic} considers some
''iterations`` of Flett's auxiliary function in terms of an
infinitesimal function, i.e., for $f\in\mathcal{D}\langle
a,b\rangle$ which is differentiable arbitrary number of times in a
right neighbourhood of the point $a$ he defines the following
functions
$$\alpha_1(x) =
\begin{cases}
\phantom{}_a^x\mathcal{K}(f)-f'(a), & x\in(a,b\rangle \\
0, & x=a
\end{cases},\quad\dots\quad
\alpha_{k+1}(x) =
\begin{cases}
\phantom{}_a^x\mathcal{K}(\alpha_k)-\alpha_k'(a), & x\in(a,b\rangle \\
0, & x=a,
\end{cases}
$$ for $k=1,2,\dots$ Then he proves the following result.

\begin{theorem}[Male\v{s}evi\`{c}, 1999]
Let $f\in\mathcal{D}\langle a,b\rangle$ and $f$ be $(n+1)$-times
differentiable in a right neighbourhood of the point $a$. If
$$\mathrm{either}\quad \alpha_n'(b) \alpha_n(b)<0, \quad \mathrm{or} \quad \alpha_n'(a) \alpha_n(b)<0,$$
then there exists $\eta\in(a,b)$ such that $\alpha_n'(\eta)=0$.
\end{theorem}

For $n=1$ Male\v{s}evi\`{c}'s condition $\alpha_1'(b)
\alpha_1(b)<0$ is equivalent to Trahan's
condition~(\ref{predpoklad}) where the second differentiability of
$f$ in a right neighbourhood of $a$ is a superfluous constraint.
The second Male\v{s}evi\`{c}'s condition $\alpha_1'(a)
\alpha_1(b)<0$ is equivalent to our condition~(\ref{malesevic}),
because $$\alpha_1'(a)=\lim_{x\to a^+}
\frac{f(x)-f(a)-f'(a)(x-a)}{(x-a)^2} = \lim_{x\to a^+}
\frac{f'(x)-f'(a)}{2(x-a)} = \frac{1}{2}f''(a)$$ and then the
inequality
$$0>\alpha_1'(a) \alpha_1(b) = \frac{1}{2}f''(a)\left(\phantom{}_a^b\mathcal{K}(f)-f'(a)\right)$$
holds if and only if~(\ref{malesevic}) holds. However, we require
only the existence of $f''(a)$ in~(\ref{malesevic}). Note that for
$n>1$ Male\v{s}evi\`{c}'s result does not correspond to
Pawlikowska's theorem (a generalization of Flett's theorem for
higher-order derivatives), see Section~\ref{secPawlikowska}, but
it goes a different way.

\begin{figure}
\begin{center}
\includegraphics[scale=0.8]{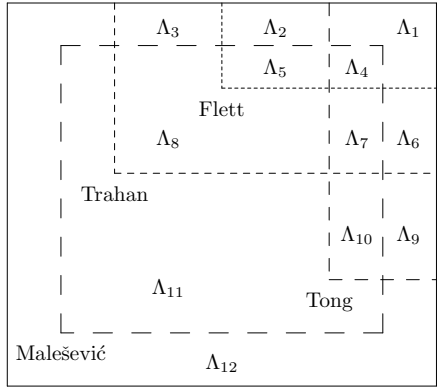}
\caption{The diagram of four families of
functions}\label{obrdiagramy2}
\end{center}
\end{figure}

Fig.~\ref{obrdiagramy2} shows all the possible cases of relations
of classes of functions satisfying assumptions of Flett, Trahan,
Tong and Male\v{s}evi\'{c}, respectively. Some examples of
functions belonging to sets $\Lambda_1,\dots,\Lambda_{12}$ were
already mentioned (e.g. $y=x^3$, $x\in\langle -1,1\rangle$,
belongs to $\Lambda_1$; $y=\sin x$,
$x\in\langle-\frac{\pi}{2},\frac{5}{2}\pi\rangle$, belongs to
$\Lambda_2$; $y=x^3$, $x\in\langle -\frac{2}{3},1\rangle$, belongs
to $\Lambda_3$; $y=\arcsin x$, $x\in\langle -1,1\rangle$, belongs
to $\Lambda_9$, and $y=\textrm{sgn}\,x$, $x\in\langle
-2,1\rangle$, belongs to $\Lambda_{12}$), other (and more
sophisticated) examples is not so difficult to find.

\begin{remark}\rm
Again, if we strengthen our assumption and consider only the
functions $f\in\mathcal{D}\langle a,b\rangle$ which are twice
differentiable at $a$, we may ask the legitimate question:
\textit{Is each of the sets $\Lambda_i$, $i=1,\dots,12$, in
Fig.~\ref{obrdiagramy2} non-empty?} In the positive case, it would
be interesting to provide a complete characterization of all the
classes of functions.
\end{remark}

\begin{problem}
All the presented conditions are only sufficient for the assertion
of Flett's theorem to hold. Provide necessary condition(s) for the
validity of~(\ref{flettMV}).
\end{problem}

\section{Integral Flett's mean value theorem}\label{secintegral}

Naturally as in the case of Lagrange's theorem we may ask whether
Flett's theorem has its analogical form in integral calculus.
Consider therefore a function $f\in\mathcal{C}\langle a,b\rangle$.
Putting
$$F(x)=\int_{a}^{x}f(t)\,\mathrm{d}t, \quad x\in\langle a,b\rangle,$$
the fundamental theorem of integral calculus yields that
$F\in\mathcal{D}\langle a,b\rangle$ with $F'(a)=f(a)$ and
$F'(b)=f(b)$. If $f(a)=f(b)$, then the function $F$ on the
interval $\langle a,b\rangle$ fulfils the assumptions of Flett's
theorem and we get the following result. It was proved by
\textsc{Stanley G. Wayment} in 1970 and it is nothing but the
integral version of Flett's theorem. Although our presented
reflection is a trivial proof of this result, we add here the
original Wayment's proof adopted from~\cite{Wayment} which does
not use the original Flett's theorem.

\begin{theorem}[Wayment, 1970]
If $f\in\mathcal{C}\langle a,b\rangle$ with $f(a)=f(b)$, then
there exists $\eta\in(a,b)$ such that {\setlength\arraycolsep{2pt}
\begin{eqnarray}\label{wayment}
f(\eta)=I_f(a,\eta).
\end{eqnarray}}
\end{theorem}

\proof Consider the function
$$F(t)=\begin{cases}
(t-a)[f(t)-I_f(a,t)], & t\in(a,b\rangle, \\ 0, & t=a.
\end{cases}$$ If $f$ is a constant on $\langle a,b\rangle$, then
$F\equiv 0$ and the assertion of theorem holds trivially. Thus,
suppose that $f$ is non-constant. Since $f\in\mathcal{C}\langle
a,b\rangle$, then by Weierstrass' theorem on the existence of
extrema there exist points $t_{1}, t_{2}\in\langle a,b\rangle$
such that
$$(\forall t\in\langle a,b\rangle)\,\, f(t_{1})\leq f(t)\leq f(t_{2}).$$ From $f(a)=f(b)$ we deduce that $f$
cannot achieve both extrema at the end points $a$ and $b$.

If $t_{2}\neq a$, then $F(b)<0<F(t_{2})$ and by Bolzano's theorem
there exists $\eta\in (t_2,b)$ such that $F(\eta)=0$. If
$t_{1}\neq a$, then $F(t_{1})<0<F(b)$ and analogously as above we
conclude that there exists $\eta\in(t_{1},b)$ such that
$F(\eta)=0$. Finally, consider the case when none of $t_{1}$ and
$t_{2}$ is equal to $a$. Then
$$a<\min\{t_{1},t_{2}\}<\max\{t_{1},t_{2}\}<b,$$ and so $F(t_{1})\leq
0\leq F(t_{2})$. From Bolzano's theorem applied to function $F$ on
the interval $\langle t_1,t_2\rangle$ we have that there exists a
point $\eta\in (a,b)$ such that $F(\eta)=0$. \qed

\begin{figure}
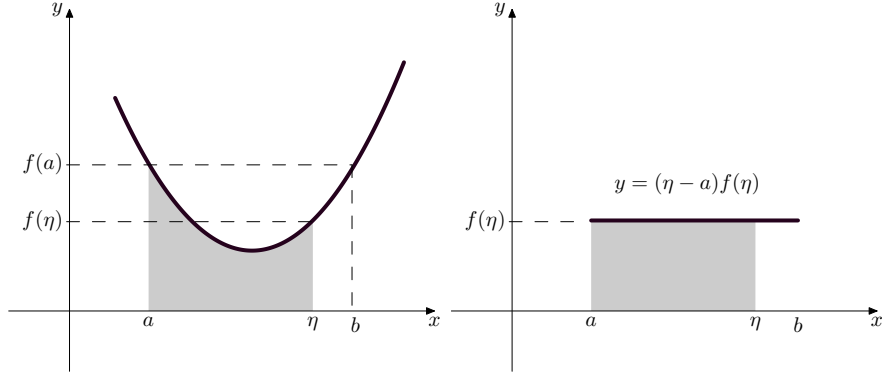

\begin{center}
\includegraphics[scale=0.8]{integralFlett.1}
\includegraphics[scale=0.8]{integralFlett2.1}
\caption{Geometrical interpretation of Wayment's theorem}
\label{obrintegralnaFlett}
\end{center}
\end{figure}

\paragraph{Geometrical meaning of Wayment's theorem}
Geometrically Wayment's theorem says that the area under the curve
$f$ on the interval $\langle a,\eta\rangle$ is equal to $(\eta-a)
f(\eta)$, i.e., volume of rectangle with sides $\eta-a$ and
$f(\eta)$, see Fig.~\ref{obrintegralnaFlett}.

Removing the condition $f(a)=f(b)$ yields the following integral
version of Riedel-Sahoo's theorem. Its proof is based on using
Riedel-Sahoo's theorem for function
$$F(x)=\int_{a}^{x}f(t)\,\mathrm{d}t, \quad x\in\langle
a,b\rangle.$$

\begin{theorem}
If $f\in\mathcal{C}\langle a,b\rangle$, then there exists
$\eta\in(a,b)$ such that
$$f(\eta)= I_f(a,\eta)
+\frac{\eta-a}{2}\cdot\phantom{}_a^b\mathcal{K}(f).
$$
\end{theorem}

In what follows we present some results from
Section~\ref{sectionsufficient} in their integral form to show
some sufficient conditions for validity of~(\ref{wayment}) with a
short idea of their proofs. The first one is Trahan's result.

\begin{proposition}
If $f\in\mathcal{C}\langle a,b\rangle$ and
$$\left[f(a)-I_f(a,b)\right]\cdot\left[f(b)-I_f(a,b)\right]\geq
0,$$ then there exists $\eta\in(a,b\rangle$ such
that~(\ref{wayment}) holds.
\end{proposition}

For the proof it is enough to consider the function
\begin{equation} g(x) =
\begin{cases}
I_f(a,x), & x\in(a,b\rangle \nonumber\\
f(a), & x=a,
\end{cases}
\end{equation}
and apply Trahan's lemma~\cite[Lemma~1]{trahan}. To show an
analogy with Tong's result consider the following means
$$B_f(a,b) = \frac{b-a}{2} I_f(a,b), \qquad
J_f(a,b) = b\,I_f(a,b)-\frac{1}{b-a}\int_{a}^{b}t
f(t)\,\mathrm{d}t.$$

\begin{proposition}\label{novytong}
Let $f\in\mathcal{C}\langle a,b\rangle$. If $B_f(a,b)=J_f(a,b)$,
then there exists $\eta\in(a,b)$ such that $f(\eta)=I_f(a,\eta)$.
\end{proposition}

In the proof we consider the auxiliary function
$$h(x)=
\begin{cases}
(x-a)[B_f(a,x)-J_f(a,x)], & x\in (a,b\rangle, \\ 0, &
x=a,\end{cases}$$ and the further steps coincide with the original
Tong's proof of Theorem~\ref{tong}.

Removing the condition $B_f(a,b)=J_f(a,b)$ we obtain the following
result.

\begin{proposition}
Let $f\in\mathcal{C}\langle a,b\rangle$. Then there exists
$\eta\in(a,b)$ such that
$$f(\eta)=I_f(a,\eta)+\frac{6[B_f(a,b)-J_f(a,b)]}{(b-a)^{2}}(\eta-a).$$
\end{proposition}

Proof, in which we use the following auxiliary function
$$H(x)=f(x)-\frac{6[B_f(a,b)-J_f(a,b)]}{(b-a)^{2}}(2x-a-b), \quad
x\in \langle a,b\rangle,$$ is again analogous to the proof of
Tong's Theorem~\ref{tong2}.

In the end of this chapter we formally present integral analogies
of new sufficient conditions of validity of Flett's theorem.

\begin{lemma}\label{lema1*}
Let $f, g\in\mathcal{C}\langle a,b\rangle$ and
\begin{equation}\label{intpodm}
\left[f(a) I_g(a,b)-g(a)I_f(a,b)\right]\cdot \left[f(b)
I_g(a,b)-g(b)I_f(a,b)\right]\geq 0.\end{equation}Then there exists
$\xi\in(a,b)$ such that
$$I_f(a,\xi)\,I_g(a,b) = I_g(a,\xi)\,I_f(a,b).$$
\end{lemma}

\proof Considering the function
$$\varphi(x)=
\begin{cases}
I_f(a,x) I_g(a,b)-I_g(a,x) I_f(a,b), & x\in (a,b\rangle, \\ 0, &
x=a,
\end{cases}$$ we have
$$\varphi'(a) =
f(a) I_g(a,b)-g(a)I_f(a,b), \quad \varphi'(b) = f(b)
I_g(a,b)-g(b)I_f(a,b).$$ Further steps are analogous to the proof
of Lemma~\ref{lema1}.\qed

\begin{theorem}
Let $f, g\in\mathcal{C}\langle a,b\rangle$ and the
inequality~(\ref{intpodm}) holds. Then there exists $\eta\in(a,b)$
such that
\begin{equation}\label{intnew}
I_g(a,b)\cdot\left(f(\eta)-I_f(a,\eta)\right) =
I_f(a,b)\cdot\left(g(\eta)-I_g(a,\eta)\right).
\end{equation}
\end{theorem}

\proof Take the auxiliary function
$$F(x) = \begin{cases}
I_g(a,b)\,I_f(a,x)-I_f(a,b)\,I_g(a,x),
 & x\in (a,b\rangle,\\
f(a) I_g(a,b)-g(a) I_f(a,b),\quad & x=a.
\end{cases}$$ By Lemma~\ref{lema1*} there exists $\xi\in(a,b)$ such that $F(\xi)=0=F(b)$. Then
by Rolle's theorem for $F$ on the interval $\langle \xi,b\rangle$
there exists $\eta\in(\xi,b)$ such that $F'(\eta)=0$.\qed

Similarly we may prove the following integral versions of
Lemma~\ref{lema2} and Theorem~\ref{veta2}.

\begin{lemma}\label{lema2*}
Let $f, g\in\mathcal{C}\langle a,b\rangle$ and $f, g$ be
differentiable at $a$. If
\begin{equation}\label{intpodm2}
\left[f(a) I_g(a,b)-g(a) I_f(a,b)\right]\cdot\left[f'(a)
I_g(a,b)-g'(a) I_f(a,b)\right]> 0,
\end{equation}
then there exists $\xi\in(a,b)$ such that
$$I_g(a,b)\cdot\left(f(a)-I_f(a,\xi)\right)
= I_f(a,b)\cdot\left(g(a)-I_g(a,\xi)\right).$$
\end{lemma}

\begin{theorem}\label{veta2*}
Let $f, g\in\mathcal{C}\langle a,b\rangle$ and $f, g$ are
differentiable at $a$. If~(\ref{intpodm2}) holds, then there
exists $\eta\in(a,b)$ such that~(\ref{intnew}) holds.
\end{theorem}

\section{Flett's theorem for higher-order
derivatives}\label{secPawlikowska}

The previous sections dealt with the question of replacing the
condition $f(a)=f(b)$ in Rolle's theorem by $f'(a)=f'(b)$. In this
section we will consider a natural question of generalizing
Flett's theorem for higher-order derivatives. We will provide the
original solution of Pawlikowska and present a new proof of her
result together with some other observations.

\subsection{Pawlikowska's theorem}

The problem of generalizing Flett's theorem to higher-order
derivatives was posed first time by \textsc{Zsolt Pales} in 1997
at the 35th international symposium of functional equations in
Graz. Solution has already appeared two years later by polish
mathematician \textsc{Iwona Paw\-likowska} in her
paper~\cite{Pawlikowska} and it has the following form.

\begin{theorem}[Pawlikowska, 1999]\label{Pawlikowska}
If $f\in\mathcal{D}^n\langle a,b\rangle$ with
$f^{(n)}(a)=f^{(n)}(b)$, then there exists $\eta\in(a,b)$ such
that
\begin{equation}\label{npawlikowska*}
\phantom{}_a^{\eta}\mathcal{K}(f) =
\sum_{i=1}^{n}\frac{(-1)^{i+1}}{i!}(\eta-a)^{i-1}f^{(i)}(\eta).
\end{equation}
\end{theorem}

Pawlikowska in her paper~\cite{Pawlikowska} generalized original
Flett's proof in such a way that she uses $(n-1)$-th derivative of
Flett's auxiliary function $g$ given by~(\ref{flettfunkcia}) and
Rolle's theorem. More precisely, the function
\begin{equation} \label{npaw2} G_{f}(x) =
\begin{cases}
g^{(n-1)}(x), & x\in(a,b\rangle \\
\frac{1}{n}f^{(n)}(a), & x=a.
\end{cases}
\end{equation} plays here an important role. Indeed, $G_{f}\in\mathcal{C}\langle
a,b\rangle\cap\mathcal{D}(a,b)$ and $$g^{(n)}(x) =
\frac{(-1)^{n}n!}{(x-a)^{n}}\left(\phantom{}_a^{x}\mathcal{K}(f)+\sum_{i=1}^{n}\frac{(-1)^{i}}{i!}(x-a)^{i-1}f^{(i)}(x)\right)
= \frac{1}{x-a}\left(f^{(n)}(x)-n\,g^{(n-1)}(x)\right)$$ for
$x\in(a,b\rangle$ which can be verified by induction. Moreover, if
$f^{(n+1)}(a)$ exists, then
$$\lim_{x\to a^+}g^{(n)}(x) = \frac{1}{n+1} f^{(n+1)}(a).$$
Further steps of Pawlikowska's proof is analogous to the original
proof of Flett's theorem using Rolle's theorem. Similarly we may
proceed using Fermat's theorem.

We have found a new proof of Pawlikowska theorem (it was not
published yet, it exists only in the form of
preprint~\cite{Molnarova}) which deals only with Flett's theorem.
The basic idea consists in iteration of Flett's theorem using an
appropriate auxiliary function.

\vskip 10pt\noindent \textbf{New proof of Pawlikowska's
theorem.}\,\,\,\, For $k=1,2,...,n$ consider the function
$$\varphi_{k}(x)=\sum_{i=0}^{k}\frac{(-1)^{i+1}}{i!}(k-i)(x-a)^{i}f^{(n-k+i)}(x)+xf^{(n-k+1)}(a), \quad x\in\langle a,b\rangle.$$
Running through all indices $k=1, 2, \dots, n$ we show that its
derivative fulfills assumptions of Flett's mean value theorem and
it implies the validity of Flett's mean value theorem for $l$-th
derivative of $f$, where $l=n-1, n-2, \dots, 1$.

Indeed, for $k=1$ we have $$\varphi_1(x)=-f^{(n-1)}(x)+x
f^{(n)}(a)\quad \textrm{and} \quad
\varphi_1'(x)=-f^{(n)}(x)+f^{(n)}(a)$$ for each $x\in\langle
a,b\rangle$. Clearly, $\varphi_1'(a)=0=\varphi_1'(b)$, so applying
Flett's theorem for $\varphi_1$ on $\langle a,b\rangle$ there
exists $u_1\in(a,b)$ such that
\begin{equation}\label{n-1}
\varphi_1'(u_1) =  \phantom{}_a^{u_1}\mathcal{K}(\varphi_1)\quad
\Leftrightarrow \quad
\phantom{}_a^{u_1}\mathcal{K}\left(f^{(n-1)}\right) =
f^{(n)}(u_1).
\end{equation} Then for $\varphi_2(x)=-2f^{(n-2)}(x)+(x-a) f^{(n-1)}(x)+x
f^{(n-1)}(a)$ we get
$$\varphi_2'(x)=-f^{(n-1)}(x)+(x-a)f^{(n)}(x)+f^{(n-1)}(a)$$ and $\varphi_2'(a)=0=\varphi_2'(u_1)$ by~(\ref{n-1}).
So, by Flett's theorem for $\varphi_2$ on $\langle a,u_1\rangle$
there exists $u_2\in (a,u_1)\subset (a,b)$ such that
$$ \varphi_2'(u_2) =  \phantom{}_a^{u_2}\mathcal{K}(\varphi_2)\quad
\Leftrightarrow \quad
\phantom{}_a^{u_2}\mathcal{K}\left(f^{(n-2)}\right) =
f^{(n-1)}(u_2) - \frac{1}{2}(u_2-a)f^{(n)}(u_2).$$ Continuing this
way after $n-1$ steps, $n\geq 2$, there exists $u_{n-1}\in (a,b)$
such that
\begin{equation}\label{1}
\phantom{}_a^{u_{n-1}}\mathcal{K}(f') = \sum_{i=1}^{n-1}
\frac{(-1)^{i+1}}{i!} (u_{n-1}-a)^{i-1} f^{(i+1)}(u_{n-1}).
\end{equation} Considering the function $\varphi_n$ we get
$$\varphi_n'(x) =
-f'(x)+f'(a)+\sum_{i=1}^{n-1}\frac{(-1)^{i+1}}{i!}(x-a)^i
f^{(i)}(x) = f'(a)+\sum_{i=0}^{n-1} \frac{(-1)^{i+1}}{i!}(x-a)^i
f^{(i+1)}(x).$$ Clearly, $\varphi_n'(a)=0$ and then
$$\phantom{}_a^{u_{n-1}}\mathcal{K}(\varphi_n') = %\frac{\varphi_n'(u_{n-1})}{u_{n-1}-a} =
-\phantom{}_a^{u_{n-1}}\mathcal{K}(f') + \sum_{i=1}^{n-1}
\frac{(-1)^{i+1}}{i!} (u_{n-1}-a)^{i-1} f^{(i+1)}(u_{n-1})=0$$
by~(\ref{1}). From it follows that $\varphi_n'(u_{n-1})=0$ and by
Flett's theorem for $\varphi_n$ on $\langle a,u_{n-1}\rangle$
there exists $\eta\in (a,u_{n-1})\subset (a,b)$ such that
\begin{equation}\label{FMV1}
\varphi_n'(\eta) = \phantom{}_a^{\eta}\mathcal{K}(\varphi).
\end{equation} Since
$$\varphi_n'(\eta) = f'(a)+\sum_{i=1}^{n}
\frac{(-1)^{i}}{(i-1)!}(\eta-a)^{i-1} f^{(i)}(\eta)$$ and
$$\phantom{}_a^{\eta}\mathcal{K}(\varphi) =
f'(a)-n\cdot\phantom{}_a^{\eta}\mathcal{K}(f)+\sum_{i=1}^{n}
\frac{(-1)^{i+1}}{i!}(n-i)(\eta-a)^{i-1} f^{(i)}(\eta),$$ the
equality~(\ref{FMV1}) yields
$$-n\cdot\phantom{}_a^{\eta}\mathcal{K}(f) = \sum_{i=1}^{n}
\frac{(-1)^{i}}{(i-1)!}(\eta-a)^{i-1}
f^{(i)}(\eta)\left(1+\frac{n-i}{i}\right) = n \sum_{i=1}^{n}
\frac{(-1)^{i}}{i!}(\eta-a)^{i-1} f^{(i)}(\eta),$$ which
corresponds to~(\ref{npawlikowska*}).\qed

\begin{remark}\rm
Recall that the assertion of Flett's theorem has an equivalent
form $f(a)=T_1(f,\eta)(a)$. Now, in the assertion of Pawlikowska's
theorem we can observe a deeper (and very natural) relation with
Taylor's polynomial. Indeed, $f(a)=T_n(f,\eta)(a)$ is an
equivalent form of~(\ref{npawlikowska*}). Geometrically it means
that Taylor's polynomial $T_n(f,\eta)(x)$ intersects the graph of
$f$ at the point $A=[a,f(a)]$.
\end{remark}

\begin{remark}\rm\label{pozndet}
Equivalent form of Pawlikowska's theorem in the form of
determinant is as follows $$\left|\begin{array}{ccccc} f^{(n)}(\eta) & h_n^{(n)}(\eta) & h_{n-1}^{(n)}(\eta) & \dots & h_0^{(n)}(\eta) \\
f^{(n-1)}(\eta) & h_n^{(n-1)}(\eta) & h_{n-1}^{(n)}(\eta) & \dots & h_0^{(n-1)}(\eta) \\
\vdots & \vdots &\vdots &\dots &\vdots \\
f^{'}(\eta) & h_n^{'}(\eta) & h_{n-1}^{'}(\eta) & \dots & h_0^{'}(\eta) \\
f(\eta) & h_n(\eta) & h_{n-1}(\eta) & \dots & h_0(\eta) \\
f(a) & h_n(a) & h_{n-1}(a) & \dots & h_0(a) \\
\end{array}\right| = 0, \qquad h_i(x)=\frac{x^i}{i!}.$$ Verification of this fact is not
as complicated as it is rather long. It is based on $n$-times
application of Laplace's formula according to the last column.
\end{remark}

Again we may ask whether it is possible to remove the condition
$f^{(n)}(a)=f^{(n)}(b)$ to obtain Lagrange's type result. The
first proof of this fact was given in Pawlikowska's
paper~\cite{Pawlikowska}. Here we present two other new proofs.

\begin{theorem}[Pawlikowska, 1999]\label{pawlikowska2}
If $f\in\mathcal{D}^n\langle a,b\rangle$, then there exists
$\eta\in (a,b)$ such that
$$f(a) = T_n(f,\eta)(a)+\frac{(a-\eta)^{n+1}}{(n+1)!}
\cdot\phantom{}_a^b\mathcal{K}\left(f^{(n)}\right).$$
\end{theorem}

\vskip 10pt\noindent \textbf{Proof I.}\,\,\,\, Consider the
auxiliary function
$$\psi_{k}(x)=\varphi_{k}(x)+\frac{(-1)^{k+1}(x-a)^{k+1}}{(k+1)!}\cdot\phantom{}_a^b\mathcal{K}\left(f^{(n)}\right), \quad k=1,2...,n,$$
where $\varphi_{k}$ is the function from our proof of
Pawlikowska's theorem. Then
$$\psi_{1}(x)=\varphi_{1}(x)+\frac{(x-a)^{2}}{2}\cdot\phantom{}_a^b\mathcal{K}\left(f^{(n)}\right),\quad x\in\langle a,b\rangle,$$
and so
$$\psi_{1}'(x)=\varphi_{1}'(x)+(x-a)\cdot\phantom{}_a^b\mathcal{K}\left(f^{(n)}\right).$$
Thus, $\psi_{1}'(a)=0=\psi_{1}'(b)$ and by Flett's theorem for
function $\psi_1$ on $\langle a,b\rangle$ there exists
$u_{1}\in(a,b)$ such that
$$\phantom{}_a^{u_1}\mathcal{K}\left(f^{(n-1)}\right)=f^{(n)}(u_{1})+\frac{a-u_{1}}{2}\cdot\phantom{}_a^b\mathcal{K}\left(f^{(n)}\right).$$
After $n-1$ steps we conclude that there exists $u_{n-1}\in(a,b)$
such that
$$
\phantom{}_a^{u_{n-1}}\mathcal{K}\left(f'\right)=
\sum_{i=1}^{n-1}\frac{(-1)^{i+1}}{i!}(u_{n-1}-a)^{i-1}f^{(i+1)}(u_{n-1})-
\frac{(a-u_{n-1})^{n-1}}{n!}\cdot\phantom{}_a^b\mathcal{K}\left(f^{(n)}\right).
$$ By Flett's theorem for $\psi_{n}$ on the interval $\langle a,u_{n-1}\rangle$
we get the desired result (all the steps are identical with the
steps of previous proof).\qed

\vskip 10pt\noindent \textbf{Proof II.}\,\,\,\, Applying
Pawlikowska's theorem to the function
$$F(x)=\left|\begin{array}{cccc} f(x) & x^{n+1} & x^n & 1 \\ f(a) & a^{n+1} & a^n & 1 \\
f^{(n)}(a) & (n+1)!\,a & n! & 0 \\ f^{(n)}(b) & (n+1)!\,b & n! & 0
\end{array}\right|,\,\, x\in\langle a,b\rangle$$ we get the result. \qed

\begin{remark}\rm
By Remark~\ref{pozndet} we may rewrite the assertion of
Theorem~\ref{pawlikowska2} in the form of determinant as follows
$$\left|\begin{array}{ccccc} f^{(n)}(\eta) & h_n^{(n)}(\eta) & h_{n-1}^{(n)}(\eta) & \dots & h_0^{(n)}(\eta) \\
f^{(n-1)}(\eta) & h_n^{(n-1)}(\eta) & h_{n-1}^{(n)}(\eta) & \dots & h_0^{(n-1)}(\eta) \\
\vdots & \vdots &\vdots &\dots &\vdots \\
f^{'}(\eta) & h_n^{'}(\eta) & h_{n-1}^{'}(\eta) & \dots & h_0^{'}(\eta) \\
f(\eta) & h_n(\eta) & h_{n-1}(\eta) & \dots & h_0(\eta) \\
f(a) & h_n(a) & h_{n-1}(a) & \dots & h_0(a) \\
\end{array}\right| = \phantom{}_a^b\mathcal{K}\left(f^{(n)}\right)\cdot\frac{(\eta-a)^n}{(n+1)!},$$ where $h_i(x)=\frac{x^i}{i!}$
for $i=0,1,\dots,n$.
\end{remark}

A relatively easy generalization of Pawlikowska's theorem may be
obtained for two functions (as a kind of "Cauchy`` version of it).

\begin{theorem}\label{pawlikowskafg}
Let $f, g\in\mathcal{D}^n\langle a,b\rangle$ and $g^{(n)}(a)\neq
g^{(n)}(b)$. Then there exists $\eta\in(a,b)$ such that
$$f(a)-T_n(f,\eta)(a)=
\phantom{}_a^b\mathcal{K}\left(f^{(n)},g^{(n)}\right)\cdot[g(a)-T_n(g,\eta)(a)].$$
\end{theorem}

\vskip 10pt\noindent \textbf{Proof I.}\,\,\,\, Considering the
function
$$h(x)=f(x)-\phantom{}_a^b\mathcal{K}\left(f^{(n)},g^{(n)}\right)\cdot g(x), \quad x\in \langle a,b\rangle,$$
we have $h\in\mathcal{D}^n\langle a,b\rangle$ and
$h^{(n)}(a)=h^{(n)}(b)$. By Pawlikowska's theorem there exists
$\eta\in(a,b)$ such that $$\phantom{}_a^{\eta}\mathcal{K}(h) =
\sum_{i=1}^{n}\frac{(-1)^{i+1}}{i!}(\eta-a)^{i-1}h^{(i)}(\eta),$$
which is equivalent to the stated result. \qed

Naturally, as in the case of Flett's theorem we would like to
generalize e.g. Trahan's result for higher-order derivatives. Our
idea of this generalization is based on application of Trahan's
approach to Pawlikowska's auxiliary function~(\ref{npaw2}).

\begin{theorem}
Let $f\in\mathcal{D}^n\langle a,b\rangle$ and
$$\left(\frac{f^{(n)}(a)(a-b)^{n}}{n!}+\mathfrak{T}_f(a)\right)
\left(\frac{f^{(n)}(b)(a-b)^{n}}{n!}+\mathfrak{T}_f(a)\right)\geq
0,$$ where $\mathfrak{T}_f(a):=T_{n-1}(f,b)(a)-f(a)$. Then there
exists $\eta\in(a,b\rangle$ such that~(\ref{npawlikowska*}) holds.
\end{theorem}

\proof Consider the function $g$ given by~(\ref{flettfunkcia}) and
function $G_{f}$ given by~(\ref{npaw2}). Clearly,
$G_{f}\in\mathcal{C}\langle a,b\rangle\cap\mathcal{D}(a,b\rangle$
and
$$g^{(n)}(x)=\frac{(-1)^n n!}{(x-a)^{n+1}}\,(T_n(f,x)(a)-f(a)), \quad x\in (a,b\rangle.$$
To apply Trahan's lemma~\cite[Lemma~1]{trahan} we need to know
signum of
$$[G_{f}(b)-G_{f}(a)]\,G_{f}'(b)=\left(g^{(n-1)}(b)-\frac{1}{n}f^{(n)}(a)\right)g^{(n)}(b),$$
i.e.,
$$
-\frac{n!(n-1)!}{(b-a)^{2n+1}}\left(\frac{f^{(n)}(a)(a-b)^{n}}{n!}+\mathfrak{T}_f(a)\right)
\cdot
\left(\frac{f^{(n)}(b)(a-b)^{n}}{n!}+\mathfrak{T}_f(a)\right) \leq
0$$ by assumption. Then by Trahan's Lemma~\cite[Lemma~1]{trahan}
there exists $\eta\in(a,b\rangle$ such that $G_f'(\eta)=0$, which
is equivalent to the assertion of theorem. \qed

Here we present another proof of "Cauchy`` type
Theorem~\ref{pawlikowskafg} which is independent on Flett's
theorem, but uses again Trahan's lemma~\cite[Lemma~1]{trahan}.

\vspace{0.3cm}\noindent\textbf{Proof of
Theorem~\ref{pawlikowskafg} II.\,\,\,} For $x\in (a,b\rangle$ put
$\varphi(x)=\phantom{}_a^x\mathcal{K}(f)$ and
$\psi(x)=\phantom{}_a^x\mathcal{K}(g)$. Consider the auxiliary
function
\begin{equation}\label{fcia} F(x) =
\begin{cases}
\varphi^{(n-1)}(x)-\phantom{}_a^b\mathcal{K}\left(f^{(n)},
g^{(n)}\right)\cdot\psi^{(n-1)}(x),
& x\in(a,b\rangle \nonumber\\
\frac{1}{n}\left[f^{(n)}(a)-\phantom{}_a^b\mathcal{K}\left(f^{(n)},
g^{(n)}\right)\cdot g^{(n)}(a)\right], & x=a.\nonumber
\end{cases}
\end{equation}
Then $F\in\mathcal{C}\langle a,b\rangle\cap\mathcal{D}(a,b\rangle$
and for $x\in(a,b\rangle$ we have {\setlength\arraycolsep{2pt}
\begin{eqnarray*}
F'(x) & = &
\varphi^{(n)}(x)-\phantom{}_a^b\mathcal{K}\left(f^{(n)},
g^{(n)}\right)\cdot\psi^{(n)}(x)
\\ & = & \frac{(-1)^n
n!}{(x-a)^{n+1}}\,\left(T_n(f,x)(a)-f(a)-\phantom{}_a^b\mathcal{K}\left(f^{(n)},
g^{(n)}\right)\cdot(T_n(g,x)(a)-g(a))\right).
\end{eqnarray*}}Then it is easy to verify that
$$[F(b)-F(a)]\,F'(b)=-\frac{1}{b-a}(F(b)-F(a))^{2}\leq 0,$$
thus by Trahan's lemma~\cite[Lemma~1]{trahan} there exists
$\eta\in(a,b\rangle$ such that $F'(\eta)=0$, i.e.,
$$\hspace{1.7cm} f(a)-T_n(f,\eta)(a)=
\phantom{}_a^b\mathcal{K}\left(f^{(n)},
g^{(n)}\right)\cdot[g(a)-T_n(g,\eta)(a)].\hspace{2cm}\square$$

\begin{remark}\rm
Similarly as in the case of Flett's and Riedel-Sahoo's points it
is possible to give the stability results for the so called $n$th
order Flett's and Riedel-Sahoo's points. These results were proved
by Pawlikowska in her paper~\cite{Pawlikowska2}, which is
recommended to interested reader. Also, some results which concern
the connection between polynomials and the set of ($n$th order)
Flett's points are proven in~\cite{Pawlikowska3}.
\end{remark}

\subsection{Flett's and Pawlikowska's theorem for divided differences}

To be able to state the general version of Flett's and
Pawlikowska's theorem in terms of divided differences, we
introduce the following necessary definitions and preliminary
results. For more details see~\cite{AIR}.

\begin{definition}\rm
\textit{The divided difference of a function $f:\langle
a,b\rangle\to\mathbb{R}$} at $n+1$ distinct points
$x_{0},\dots,x_{n}$ of the interval $\langle a,b\rangle$ is
defined as follows {\setlength\arraycolsep{2pt}
\begin{eqnarray*}
[x_0;f] & := & f(x_0), \\
{[x_0,x_1;f]} & := &
\phantom{}_{x_0}^{x_1}\mathcal{K}(f), \\
{[x_0,x_1,\dots,x_n;f]} & := &
\frac{[x_0,x_1,\dots,x_{n-1};f]-[x_1,x_2,\dots,x_n;f]}{x_0-x_n},
\quad n\geq 2.
\end{eqnarray*}}\end{definition} If the points $x_0,\dots,x_n$ are not distinct,
then the divided difference is defined by a limit process
$$[\underbrace{x_{0},\dots,x_{0}}_{k+1},x_{k+1},\dots,x_{n};f] :=
\lim_{x_{1},\dots,x_{k}\to x_{0}}[x_{0},x_{1},\dots,x_{n};f]$$
provided the limit exists. In particular
$$[\underbrace{c,\dots,c}_{n+1};f] :=
\lim_{x_{1},\dots,x_{n}\to c}[c,x_{1},\dots,x_{n};f].$$ %In what
%follows, all points $x_0, \dots, x_n$ in the symbol $[x_0,\dots,
%x_n;f]$ will be assumed to be pairwise distinct unless specified otherwise.
The following result plays a key role in extension of
Flett's and Pawlikowska's theorem. For its proof and more details
we refer to the paper~\cite{AIR} and references given therein.

\begin{proposition}\label{poddif5}
Let $f\in\mathcal{C}\langle a,b\rangle$ and $f$ be $n$-times
differentiable at $a$ and $b$ with $f^{(n)}(a)=f^{(n)}(b)$. Then
there exists $\eta\in(a,b)$ such that in any neighborhood of the
point $\eta$ there exist equidistant points $\eta_0<\dots<\eta_n$,
$\eta_0<\eta<\eta_n$ such that $[a,\eta_0,\dots,\eta_n;f]=0$.
\end{proposition}

Immediately, the generalized Pawlikowska's theorem for divided
differences has the following form.

\begin{theorem}[Abel-Ivan-Riedel, 2004]\label{thmPaw}
If $f\in\mathcal{D}^{n}\langle a,b\rangle$, then there exists
$\eta\in(a,b)$ such that in any neighborhood of the point $\eta$
there exist equidistant points $\eta_0<\dots<\eta_n$,
$\eta_0<\eta<\eta_n$ such that
$$[a,\eta_0,\dots,\eta_n;f]=\frac{1}{(n+1)!}\phantom{}_a^b\mathcal{K}\left(f^{(n)}\right).$$
\end{theorem}

\proof Using the relation
$\left[a,\eta_0,\dots,\eta_n;(x-a)^{n+1}\right]_x=1$ and applying
Proposition~\ref{poddif5} for the function $h:\langle
a,b\rangle\to \mathbb{R}$ given by
$$h(x)=f(x)-\frac{1}{(n+1)!}\phantom{}_a^b\mathcal{K}\left(f^{(n)}\right)(x-a)^{n+1}$$
yields the desired result. \qed

If in Theorem~\ref{thmPaw} we take $\eta_i\rightarrow \eta$ for
$i=0,\dots,n$, then we get a new form of Pawlikowska's theorem
without boundary assumption.

\begin{corollary}\rm If
$f\in\mathcal{D}^{n}\langle a,b\rangle$, then there exists
$\eta\in(a,b)$ such that
$$[a,\underbrace{\eta,\dots,\eta}_{n+1};f]=\frac{1}{(n+1)!}\phantom{}_a^b\mathcal{K}\left(f^{(n)}\right).$$
\end{corollary}

For $n=1$ this implies a new form of Flett's mean value theorem.

\begin{corollary}\rm If $f\in\mathcal{D}\langle
a,b\rangle$ and $f'(a)=f'(b)$, then there exists $\eta\in(a,b)$
such that $[a,\eta,\eta;f]=0$.
\end{corollary}

Finally, a Pawlikowska's type theorem with boundary has the
following form.

\begin{theorem}[Abel-Ivan-Riedel, 2004]
If $f\in\mathcal{D}^{n}\langle a,b\rangle$ and
$f^{(n)}(a)=f^{(n)}(b)$, then there exists $\eta\in(a,b)$ such
that
$$[a,\underbrace{\eta,\dots,\eta}_{n};f]=\frac{f^{(n)}(\eta)}{n!}.$$
\end{theorem}

\proof By Proposition~\ref{poddif5} there exists $\eta\in(a,b)$
such that in any neighbourhood of $\eta$ there exist equidistant
points $\eta_0<\dots<\eta_n$, $\eta_0<\eta<\eta_n$ such that
$[a,\eta_0,\dots,\eta_n;f]=0$. This yields
$$[a,\eta_1,\dots,\eta_n;f]-[\eta_0,\dots,\eta_n;f]=0,$$ and thus
for $\eta_i\to \eta$ for $i=0,\dots,n$ we get
$$[a,\underbrace{\eta,\dots,\eta}_n;f]=[\underbrace{\eta,\dots,\eta}_{n+1};f]
=\frac{f^{(n)}(\eta)}{n!},$$ where the last equality follows from
Stieltjes' theorem, see~\cite{AIR}.\qed

%\begin{remark}\rm
%Note that the last result has the following geometrical
%interpretation: Taylor's polynomial $T_n(f,\eta)(x)$ intersects
%the graph of $f$ at the point $A=[a,f(a)]$.
%\end{remark}

\section{Concluding remarks}

In this paper we provided a summary of results related to Flett's
mean value theorem of differential and integral calculus of a
real-valued function of one real variable. Indeed, we showed that
for $f\in\mathcal{D}\langle a,b\rangle$ the assertion of Flett's
theorem holds in each of the following cases:
\begin{itemize}
\item[(i)] $f'(a)=f'(b)$ (Flett's condition); \item[(ii)]
$(f'(a)-\phantom{}_a^b\mathcal{K}(f))\cdot(f'(b)-\phantom{}_a^b\mathcal{K}(f))\geq
0$ (Trahan's condition); \item[(iii)] $A_f(a,b)=I_f(a,b)$ (Tong's
condition); \item[(iv)] $(f'(a)-\phantom{}_a^b\mathcal{K}(f))\cdot
f''(a)>0$ provided $f''(a)$ exists (Male\v{s}evi\`{c}'s
condition).
\end{itemize}
Then we discussed possible generalization of Flett's theorem to
higher-order derivatives and provided a new proof of Pawlikowska's
theorem and related results. Up to a few questions and open
problems explicitly formulated in this paper, there are several
problems and directions for the future research.

The survey of results related to Flett's mean value theorem should
be continued in~\cite{HM}, because we did not mention here any
known and/or new generalizations and extensions of Flett's theorem
made at least in two directions: to move from the real line to
more general spaces (e.g. vector-valued functions of vector
argument~\cite{RS}, holomorphic functions~\cite{DPRS}, etc.),
and/or to consider other types of differentiability of considered
functions (e.g. Dini's derivatives~\cite{Reich}, symmetric
derivatives~\cite{Sahoo}, $v$-derivatives~\cite{Ohriska}, etc.).
Also, a characterization of all the functions that attain their
Flett's mean value at a particular point between the endpoints of
the interval~\cite{RSa}, other functional equations and means
related to Flett's theorem should be mentioned in the future.

\section*{Acknowledgement}

This research was partially supported by Grant VEGA 1/0090/13.

%%%%%%%%%%%%%%%%%%%%%%%%%%%%%%%%%%%%%%%%%%%%%%%%%%%%%%%%%%%%%%%%%%%%%%%%%%

\end{document}